\documentclass[a4paper, 12pt]{article}
\usepackage[utf8]{inputenc}
\usepackage{ marvosym }
\usepackage{cmap}
\usepackage[T2A]{fontenc}
\usepackage{amsmath,amssymb,amsthm}
\usepackage[a4paper,hmargin=2.5cm,vmargin=2.5cm]{geometry}
\usepackage[english,russian]{babel}
\usepackage{circuitikz}

\definecolor{green}{RGB}{10, 150, 0}
\usepackage[pdftex,colorlinks=true,linkcolor=blue,urlcolor=red,unicode=true,hyperfootnotes=false,bookmarksnumbered]{hyperref}
\usepackage{indentfirst}
\usepackage{graphicx}
\usepackage{amssymb}
\usepackage{amsthm}
\usepackage{amsmath}
\usepackage{cmap}
\usepackage{ dsfont }
\usepackage{enumitem}
\usepackage{color}
\usepackage[ruled,vlined,linesnumbered]{algorithm2e}
\usepackage{geometry}
\usepackage{pgf,tikz}\usetikzlibrary{arrows}

\newtheorem{theorem}{Theorem}[section]
\newtheorem{proposition}[theorem]{Proposition}
\newtheorem{lemma}[theorem]{Lemma}

\newtheorem{corollary}[theorem]{Corollary}

\theoremstyle{definition}
\newtheorem{definition}[theorem]{Definition}
\newtheorem{example}[theorem]{Example}
\newtheorem{problem}[theorem]{Problem}

\theoremstyle{remark}
\newtheorem{remark}[theorem]{Remark}

\newcommand{\rvline}{\hspace*{-\arraycolsep}\vline\hspace*{-\arraycolsep}}

\renewcommand{\ge}{\geqslant}
\renewcommand{\le}{\leqslant}

\newcommand{\notdiplom}[0]{}
\long\def\notformoebius#1\endnotformoebius{#1}

\long\def\comment#1\endcomment{}

\newcommand{\bluenew}[1]{#1}
\newcommand{\bluevarnew}[1]{#1}
\newcommand{\newremove}{}
\newcommand{\bluer}[1]{#1}
\newcommand{\bluevar}[1]{#1}

\usepackage{xpatch}
\xpretocmd{\abstract}{\selectlanguage{english}}{}{}  
\xpretocmd{\proof}{\selectlanguage{english}}{}{}
\xpretocmd{\figure}{\selectlanguage{english}}{}{}
\xpretocmd{\thebibliography}{\selectlanguage{english}}{}{}

\begin{document}

	\notdiplom{\title{Superport networks}
	\author{\notformoebius P.~Pylyavskyy, \endnotformoebius S.~Shirokovskikh\notformoebius, M.~Skopenkov\endnotformoebius}
	\date{}
	\maketitle}
        \begin{abstract}
            We study multiport networks, common in electrical engineering. They have boundary conditions different from electrical networks: the boundary vertices are split into pairs and the sum of the incoming currents is set to be zero in each pair. If one sets the voltage difference for each pair, then the incoming currents are uniquely determined. We generalize Kirchhoff's matrix-tree theorem to this setup. Different forests now contribute with different signs, making the proof subtle. In particular, we use the formula for the response matrix minors by R. Kenyon--D. Wilson, determinantal identities, and combinatorial bijections. We introduce superport networks, generalizing both ordinary networks and multiport ones.	


        \end{abstract}
	\notdiplom{\footnotetext{This work is supported by NSF grant DMS-1949896, the Ministry of Science and Higher Education of the Russian Federation (agreement no. 075–15–2022–287), and Center of Excellence for Generative AI at King Abdullah University of Science and Technology (KAUST). 
    }}

	\section{Introduction}

Mathematical theory of electrical networks is rich and well-developed. In a nutshell, it studies weighted graphs, where we apply voltages to some vertices and read off the currents flowing from them. See the introductory reviews \bluenew{\cite{Doyle-Snell-84,PS,Skopenkov-23,ZSU}}.

However, the electrical networks we face in our everyday life are different. The vertices often come in pairs so that the currents within each pair are opposite. For instance, we charge a laptop using a plug with \emph{two} pins usually, and all the current entering through one pin exits through the other one. When we plug a device into the laptop, we use a cable with \emph{two} (and often more) wires, and the currents are again opposite because there is no other way the current entering the device can leave it. 

\begin{figure}[h]
\center
 \begin{circuitikz}
           \draw [latexslim-latexslim, color=red] (0.3,0) -- (2.7,0);
           \draw [latexslim-latexslim, color=blue] (0.3,4) -- (2.7,4);
             \draw[color=red]
             (0.4,0) to[short, l_=$\Delta U_{34}$] (2.6,0);
             \draw[color=blue] (0.4,4) to[short, l=$\Delta U_{12}$] (2.6,4);
             \draw[color=blue] (-0.4, 4.4) node[flowarrow, rotate=-45]{\rotatebox{45}{}}
             (3.4, 4.4) node[flowarrow, rotate=-135]{\rotatebox{135}{$I_2$}};
             \draw[color=red] (-0.4, -0.4) node[flowarrow, rotate=45]{\rotatebox{-45}{}}
             (3.4, -0.4) node[flowarrow, rotate=135]{\rotatebox{-135}{}};
             \draw
             (0, 0) to[short, f=$I_{35}$] (0.5, 2)
             (0.5, 2) to[short, f_=$I_{56}$] (2.5, 2)
             (0, 4) to[short, f_=$I_{15}$] (0.5, 2)
             (0, 4) to[short, f=$I_{16}$] (2.5, 2)
            (3, 0) to[short, f=$I_{46}$] (2.5, 2)
            (3, 4) to[short, f=$I_{26}$] (2.5, 2)
             (0,0) node[circ, color=red]{}
             (0,4) node[circ, color=blue]{}
             (3,0) node[circ, color=red]{}
             (3,4) node[circ, color=blue]{}
            (0.5,2) node[circ]{}
             (2.5,2) node[circ]{}
            {[anchor=south west]
                            (0,0) node[color=red] {$U_3$}
                            (0,4) node[color=blue] {$U_1$}
                            (0.5,2) node[] {$U_5$}
                            (2.5,2) node[] {$U_6$}
                            }
            {[anchor=south east]
                           (3,0) node[color=red] {$U_4$}
                            (3,4) node[color=blue] {$U_2$}
                            };
                            \draw (3.7, -0.3) node[red]{$I_4$};
                            \draw (-0.7, -0.3) node[red]{$I_3$};
                            \draw (-0.6, 4.2) node[blue]{$I_1$};
             \end{circuitikz}
             \caption{Example of a multiport network. The two ports are shown in color. Given the voltage differences $\Delta U_{12}$ and $\Delta U_{34}$, the voltages $U_k$ and the currents $I_{kl}$ for $k,l=1,\dots,6$ are found from the Kirchhoff \bluevar{and Ohm} laws and the port condition $I_1=-I_2$, $I_3=-I_4$. See Definition~\ref{def-multiport}.}
             \label{fig:example}
             \end{figure}
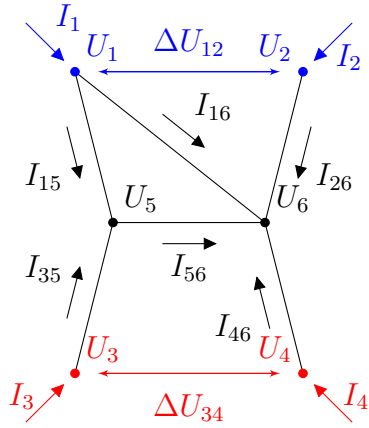

Those are examples of multiport networks. 
In electrical engineering, they are even more common than ordinary electrical networks \cite{BIK, Bess}. However, there are apparently no mathematical papers about them until now. Multiport networks have boundary conditions different from electrical networks: the boundary vertices are split into pairs and the sum of the incoming currents is set to be zero in each pair (see Fig.~\ref{fig:example}). If one sets the voltage difference between the vertices of each pair, then the incoming currents are uniquely determined. The map taking the voltage differences to the resulting currents is called the multiport network response; it describes the entire network response to external action.

The mathematical theory of electrical networks started in the middle of the 19th century with Kirchhoff's matrix-tree theorem, which is a combinatorial formula for \bluenew{the currents, hence} the response, of an electrical network (not a multiport one).
New variations and generalizations continue to appear \cite{CK}. For example, a formula for the response matrix minors was found only in 2009 by R.W.~Kenyon and D.B.~Wilson \cite{KW09}; see Theorem~\ref{Kenyon-Wilson} below. The responses of plane electrical networks are characterized by the positivity of such minors~
\cite{CM}.

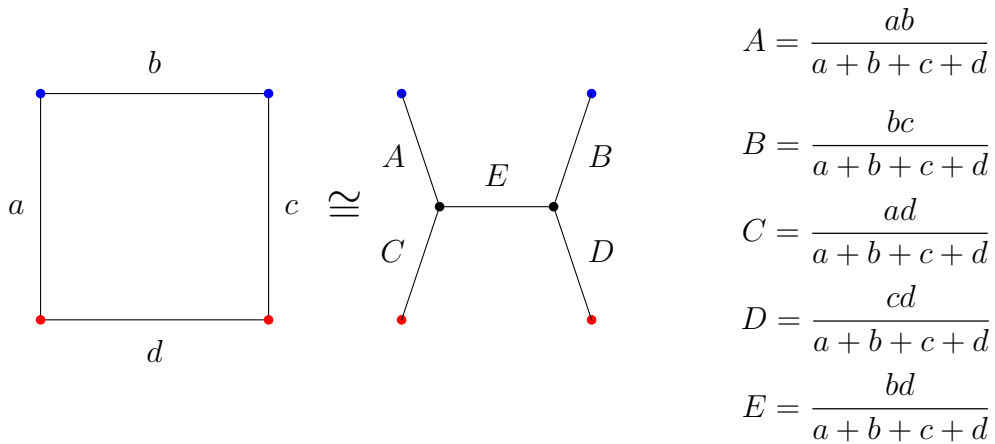
\begin{figure}[h]
\center
\begin{tabular}{cl}  
         \begin{tabular}{c}
\begin{circuitikz}
           \draw (4,1.5) node[black]{\Large$\cong$ \normalsize};
                 \draw
                 (0,0) to[short, l=$a$] (0,3)
                to[short, l=$b$] (3,3)
                to[short, l=$c$] (3,0)
                to[short, l=$d$] (0,0)
                 (0,0) node[circ, color=red]{}
                 (0,3) node[circ, color=blue]{}
                 (3,0) node[circ, color=red]{}
                 (3,3) node[circ, color=blue]{}
                 (5.25,1.5) node[circ]{}
                 (6.75,1.5) node[circ]{}
                 (7.25,3) node[circ, color=blue]{}
                 (4.75,3) node[circ, color=blue]{}
                 (4.75,0) node[circ, color=red]{}
                 (7.25,0) node[circ, color=red]{};
                 \draw (5.25,1.5) to[short, l=$E$] (6.75,1.5)
                  to[short, l_=$B$] (7.25,3);
                  \draw(5.25,1.5) to[short, l=$A$] (4.75,3);
                  \draw(5.25,1.5) to[short, l_=$C$] (4.75,0);
                  \draw(6.75,1.5) to[short, l=$D$] (7.25,0);
            \end{circuitikz}
\end{tabular}
           & \begin{tabular}{l}
             \parbox{0.3\linewidth}{             
    	$$A = \frac{ab}{a+b+c+d}$$
		$$B = \frac{bc}{a+b+c+d}$$
            $$C = \frac{ad}{a+b+c+d}$$
            $$D = \frac{cd}{a+b+c+d}$$
            $$E = \frac{bd}{a+b+c+d}$$
    }
         \end{tabular}  \\
\end{tabular}
 \caption{The Box-H transformation. It preserves the response of the two-port network and also the \emph{voltage drop} between the ports (shown in color) \cite{S}.}
	\label{fig:Box-H_transformation}
\end{figure}

Another important tool is \emph{network transformations} \cite{S}. Multiport networks have their own transformations: one example was introduced in the classical work~\cite{B} from 1927 and a new \emph{Box-H transformation} (Fig.~\ref{fig:Box-H_transformation}) was discovered by the \notformoebius\notdiplom{second }\endnotformoebius author in 2021 \cite{S}.

We introduce a new notion of \emph{superport networks}, generalizing both ordinary electrical networks and multiport networks. For this object, we prove the basic existence and uniqueness theorem of voltages, similar to the one for electrical networks (see, e.g. \cite{ZSU}) and two-port networks \cite{S}. Our main result is a generalization of Kirchhoff's matrix-tree theorem to such superport (and hence multiport) networks. See Theorems~\ref{th-Lij} and~\ref{main}. \notdiplom{It generalizes both the formula for elements of the response matrix and the one for the determinant.} In contrast to the classical Kirchhoff theorem, different forests can now contribute with different signs.
This makes the proofs subtle and requires new ideas. In particular, we use the above-mentioned formula by Kenyon \bluevarnew{and} Wilson, determinantal identities, and combinatorial bijections.

Our generalization comes with straightforward combinatorial corollaries on counting of spanning forests; see Corollary~\ref{generalized_Cayley}. 
Yet another application of superport networks is a construction of single-valued conjugate discrete harmonic functions in multiply-connected domains; see Remark~\ref{another_motivation}. All these results show that superport networks form a natural general setup to develop network theory.

This paper is organized as follows. In Section~\ref{prelim}, we recall the definition and properties of ordinary electrical networks to be used and generalized later. In Section~\ref{sec-multiport}, we introduce our main notion in the simplest nontrivial particular case. In Section~\ref{sec-exist-unique}, we define the notions in full generality and prove that they are well-defined. Section~\ref{sec-exist-unique} is independent \bluevarnew{of} the previous ones, and the reader can proceed to it immediately. In Sections~\ref{sec-response}--\ref{sec-voltages-currents}, we state and prove basic properties of superport networks. In Section~\ref{sec-matrix-tree}, we state and prove our main result (Theorem~\ref{main}). The statement of the theorem relies on Definitions~\ref{def-superport-network}, \ref{def-response}, and \ref{def-valid} only. \notformoebius In Section~\ref{sec-open} we state open problems. \endnotformoebius

 
\section{Preliminaries: electrical networks} \label{prelim}

We start by recalling the definition and properties of electrical networks. 

The following definition is borrowed almost entirely from~\cite{PS}.

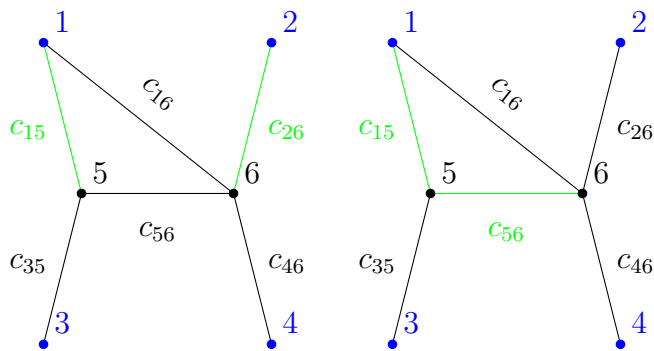
\begin{figure}[h]
 \center
  \begin{tabular}{c}
 \begin{circuitikz}
 \draw [color = green]
 (0.5, 2) to[short, l=$c_{15}$] (0, 4)
 (2.5, 2) to[short, l_=$c_{26}$] (3, 4);
             \draw
             (0, 0) to[short, l=$c_{35}$] (0.5, 2)
             (0.5, 2) to[short, l_=$c_{56}$] (2.5, 2)
             (2.5, 2) to[short, l=$c_{46}$] (3, 0)
             (2.5, 2) to[short, l_=$c_{16}$] (0, 4)
             (0,0) node[circ, color=blue]{}
             (0,4) node[circ, color=blue]{}
             (3,0) node[circ, color=blue]{}
             (3,4) node[circ, color=blue]{}
             (0.5,2) node[circ]{}
             (2.5,2) node[circ]{}
            {[anchor=south west]
                            (0,0) node[color=blue] {3}
                           (3,0) node[color=blue] {4}
                            (0,4) node[color=blue] {1}
                            (3,4) node[color=blue] {2}
                            (0.5,2) node[] {5}
                            (2.5,2) node[] {6}
                            };
            \draw[green, thick] (0,0) circle (0.15);
            \draw[green, thick] (3,0) circle (0.15);
             \end{circuitikz}
\end{tabular} \begin{tabular}{l}
            \begin{circuitikz}
 \draw [color = green]
 (0.5, 2) to[short, l_=$c_{56}$] (2.5, 2)
 (0.5, 2) to[short, l=$c_{15}$] (0, 4);
             \draw
             (0, 0) to[short, l=$c_{35}$] (0.5, 2)
             (2.5, 2) to[short, l_=$c_{16}$] (0, 4)
             (2.5, 2) to[short, l=$c_{46}$] (3, 0)
             (2.5, 2) to[short, l_=$c_{26}$] (3, 4)
             (0,0) node[circ, color=blue]{}
             (0,4) node[circ, color=blue]{}
             (3,0) node[circ, color=blue]{}
             (3,4) node[circ, color=blue]{}
             (0.5,2) node[circ]{}
             (2.5,2) node[circ]{}
            {[anchor=south west]
                            (0,0) node[color=blue] {3}
                           (3,0) node[color=blue] {4}
                            (0,4) node[color=blue] {1}
                            (3,4) node[color=blue] {2}
                            (0.5,2) node[] {5}
                            (2.5,2) node[] {6}
                            };
    \draw[green, thick] (0,0) circle (0.15);
\draw[green, thick] (3,0) circle (0.15);
\draw[green, thick] (3,4) circle (0.15);
             \end{circuitikz}
         \end{tabular}
\caption{An electrical network. The boundary vertices are \bluevarnew{1, 2, 3, 4 (in blue)}.
Examples of valid forests in the network are shown in green.}
\label{fig:valid}
             \end{figure}

\begin{definition} (See Fig.~\ref{fig:valid})
An \textit{electrical network with ${m}$ terminals} is a \bluenew{finite} connected graph with a positive number (\textit{conductance}) assigned to each edge, and \bluenew{${m \ge 1}$} marked (\textit{boundary}) vertices.
For notational convenience, we assume that the graph has neither multiple edges nor loops.

Fix an enumeration of the vertices ${1}$, ${2}$, ${\dots}$, ${n}$ such that ${1, \dots, m}$ are the boundary ones. Denote by ${kl}$ the edge between the vertices ${k}$ and ${l}$. Denote by $c_{kl}$ the conductance of the edge ${kl}$. Set $c_{kl} = 0$, if there is no edge between ${k}$ and ${l}$.

An \textit{electrical circuit} is an electrical network along with ${m}$ real numbers ${U_1, \dots, U_m}$ (\textit{\bluenew{prescribed} voltages}) assigned to the boundary vertices.
	
Each electrical circuit gives rise to certain numbers $U_k$, where $m+1\le k\le n$ (\textit{voltages} at the vertices), and $I_{kl}$, where $1\le k,l\le n$ (\textit{currents} through the edges). These numbers are determined by the following $2$ axioms:

\begin{enumerate}
     \item[(C)] \textit{The Ohm law.} For each \bluenew{ordered pair of vertices $k, l$}, we have $I_{kl} = c_{kl}(U_k - U_l)$.

	\item[(I)]\textit{The Kirchhoff current law.} For each vertex $k>m$ we have $\sum_{l=1}^{n} I_{kl}=0$.
\end{enumerate}
\end{definition}	

We emphasize that the currents $I_{kl}$ are indexed by ordered pairs of vertices ${(k, l)}$ rather than edges. In particular, by axiom~\textup{(C)} we have ${I_{kl}} = 0$ if there is no edge $kl$, and $I_{kl} = -I_{lk}.$ \bluenew{
In particular, $I_{kk} = 0$. Although there is no current flowing from a vertex to itself (and no edge), it is convenient to set $I_{kk} = 0$ for the sake of notation.}

 The numbers $U_k$ and $I_{kl}$ are well-defined by these axioms by the following classical result.

\begin{theorem} \label{prelim_th} \textup{(}See, e.g., \textup{\cite[Theorem~2.1]{PS})} For any electrical circuit
	the system of linear equations \textup{(C), (I)} in the variables $U_k$, where \bluenew{$m+1 \le k\le n$}, and $I_{kl}$, where $1\le k,l\le n$, has a unique solution.
\end{theorem}
	
The numbers $(I_{1}, \dots, I_{m}) := ({{\sum_{\bluevarnew{l}=1}^{n} I_{1\bluevarnew{l}}}}, \dots, {{\sum_{\bluevarnew{l}=1}^{n} I_{m\bluevarnew{l}}}})$ are called \textit{the incoming currents} \bluenew{(the sum of the currents through all the edges starting at a boundary vertex $k$ equals the current incoming into the network through the vertex $k$; see Fig.~\ref{fig:example})}.
The linear map $$C\colon\mathds{R}^m \to \mathds{R}^m, \quad({U_{1}}, \dots, {U_{m}}) \mapsto (I_{1}, \dots, I_{m})$$ is called \emph{the response} of the network. Denote by $C$ the matrix of this map as well, and by $C_i^j$ the entry in $i$-th row and $j$-th column. It is well-known that $C_i^j = C_j^i$. Denote by $\widetilde{C}$ the matrix obtained from $C$ by removing the last row and the last column.

 Let us state the well-known Kirchhoff's matrix-tree theorem, which allows us to compute the response \bluenew{graphically}. 

\bluenew{A spanning forest in an electrical network is \emph{valid} if each component contains exactly one boundary vertex (see Fig.~\ref{fig:valid}). In other words, we introduce the equivalence relation on the set of vertices such that all boundary  vertices are equivalent and the equivalence class of any interior vertex consists of a single vertex.} 
A spanning forest in an electrical network is valid, if it becomes a spanning tree in the quotient by this equivalence relation.
         
The \emph{weight} $w(H)$ of a forest $H$ is the product of the conductances of all the edges in the forest. Hereafter an empty product \bluenew{(respectively, sum)} is set to be $1$ \bluevarnew{(respectively, $0$)} by definition.

\begin{figure}
    \centering

\begin{tabular}{cl}

         \begin{tabular}{c}
           \begin{circuitikz}
           \draw [color = green] (0,0) to[short, l=$a$] (0,2)
                (0,2) to[short, l=$b$] (2,2)
               (2,2) to[short, l=$c$] (2,0);
                 \draw[color = lightgray]
                 (2,0) to[short, l=$d$] (0,0)
                 (0,0) node[circ, color=blue]{}
                 (0,2) node[circ, color=blue]{}
                 (2,0) node[circ, color=blue]{}
                 (2,2) node[circ, color=blue]{}
                
                {[anchor=south east]
                (2,0) node[color=blue] {4}
                (2,2) node[color=blue] {2}
                }
                {[anchor=south west]
                (0,0) node[color=blue] {3}
                (0,2) node[color=blue] {1}
                };
            \end{circuitikz}
         \end{tabular}
           & \begin{tabular}{l}
             \begin{circuitikz}
           \draw [color = green] (0,0) to[short, l=$a$] (0,2)
                (0,2) to[short, l=$b$] (2,2)
                (2,0) to[short, l=$d$] (0,0);
                 \draw[color = lightgray]
                 (2,2) to[short, l=$c$] (2,0)
                 (0,0) node[circ, color=blue]{}
                 (0,2) node[circ, color=blue]{}
                 (2,0) node[circ, color=blue]{}
                 (2,2) node[circ, color=blue]{}
                
                {[anchor=south east]
                (2,0) node[color=blue] {4}
                (2,2) node[color=blue] {2}
                }
                {[anchor=south west]
                (0,0) node[color=blue] {3}
                (0,2) node[color=blue] {1}
                };
            \end{circuitikz}  
         \end{tabular}
\end{tabular}
\begin{tabular}{cl}  
         \begin{tabular}{c}
           \begin{circuitikz}
           \draw [color = green] (0,0) to[short, l=$a$] (0,2)
                (2,0) to[short, l=$d$] (0,0)
               (2,2) to[short, l=$c$] (2,0);
                 \draw[color = lightgray]
               (0,2) to[short, l=$b$] (2,2)
                 (0,0) node[circ, color=blue]{}
                 (0,2) node[circ, color=blue]{}
                 (2,0) node[circ, color=blue]{}
                 (2,2) node[circ, color=blue]{}
                
                {[anchor=south east]
                (2,0) node[color=blue] {4}
                (2,2) node[color=blue] {2}
                }
                {[anchor=south west]
                (0,0) node[color=blue] {3}
                (0,2) node[color=blue] {1}
                };
            \end{circuitikz}
         \end{tabular}
           & \begin{tabular}{l}
             \begin{circuitikz}
           \draw [color = green] 
                (0,2) to[short, l=$b$] (2,2)
                (2,0) to[short, l=$d$] (0,0)
               (2,2) to[short, l=$c$] (2,0);
                 \draw[color = lightgray]
                 (0,0) to[short, l=$a$] (0,2)
                 (0,0) node[circ, color=blue]{}
                 (0,2) node[circ, color=blue]{}
                 (2,0) node[circ, color=blue]{}
                 (2,2) node[circ, color=blue]{}
                
                {[anchor=south east]
                (2,0) node[color=blue] {4}
                (2,2) node[color=blue] {2}
                }
                {[anchor=south west]
                (0,0) node[color=blue] {3}
                (0,2) node[color=blue] {1}
                };
            \end{circuitikz}  
         \end{tabular}
     \end{tabular}    
         $C = \left( \begin{matrix}
			a+b & -b & -a & 0\\
            -b & b+c & 0 & -c\\
            -a & 0 & a+d & -d\\
            0 & -c & -d & c+d\\
		\end{matrix}\right ) \quad \quad \det\widetilde{C} = \begingroup \color{green}abc+abd+acd+bcd\endgroup$\\
    \caption{Kirchhoff's matrix-tree theorem. The boundary vertices are \bluevarnew{1, 2, 3, 4 (in blue)}. Altogether, they form the unique valid forest. The spanning trees are in green. See Theorem~\ref{Kirchhoff}.}
    \label{fig:kirchhoff-matrix-tree}
\end{figure}
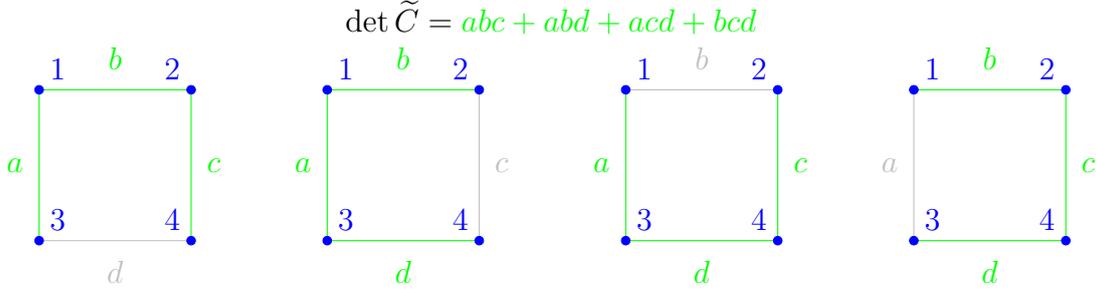

    \begin{theorem}[Kirchhoff's matrix-tree theorem] \label{Kirchhoff} \textup{(See Fig.~\ref{fig:kirchhoff-matrix-tree} \bluenew{and e.g.~\cite[Theorem~8.1]{KW09})}.)}
 	For each electrical network with $m \ge 2$ terminals \bluenew{we have:}
		\begin{enumerate}
			\item[(1)] $$\det \widetilde{C} = \frac{\sum_T w(T)}{\sum_H w(H)},$$
			where the sum in the denominator is over all valid forests $H$ in the electrical network, and the sum in the numerator is over all spanning trees $T$ in the electrical network.
			\item[(2)] \bluenew{for each $i \neq j$ we have}
            $$C_i^j = -\frac{\sum_G w(G)}{\sum_H w(H)},$$
			where the sum in the denominator is over all valid forests $H$ in the electrical network, and the sum in the numerator is over all spanning forests $G$ with $m-1$ components, such that each component contains a boundary vertex and the \bluenew{boundary} vertices $i$ and $j$ are in the same component.
		\end{enumerate}
    \end{theorem}

\bluenew{Kirchhoff~\cite{Kirchhoff-47} was originally interested in finding the currents in electrical networks. What is usually understood by Kirchhoff's matrix-tree theorem now~\cite{CK} is a closely related result, namely, a particular case of assertion (1) when all the vertices are boundary vertices (and hence the denominator equals $1$ there). In this work, more general assertion (1) and even closer-to-the-original assertion (2) are referred to as \emph{Kirchhoff's matrix-tree theorem}. Such formulation makes the dependence on the boundary conditions explicit and paves the way for generalizations considered in the present paper.}
    
Let us illustrate the theorem by an immediate combinatorial application.
    \begin{corollary}[Cayley's formula]\label{Cayley}
     The number of trees on $m$ labeled vertices is $m^{m-2}$.
    \end{corollary}	
    \begin{proof}[Sketch of the proof]
    \bluenew{Assume $m \ge 2$; otherwise, there is nothing to prove.} Consider the electrical network with $m$ terminals on the complete graph on $m$ vertices with the edges of unit conductances. Then by definition \bluenew{of the response matrix, we get}
       	 $$C_{i}^{j} = \begin{cases}
		m-1, &\text{if $i = j$},\\
		-1, &\text{if $i \neq j$},\\
	\end{cases}$$ and it is easy to find $\det \widetilde{C} = m^{m-2}.$ The only valid forest has $m$ isolated vertices, hence by Theorem~\ref{Kirchhoff} the number of spanning trees $T$ is $\sum_T 1=\sum_T w(T) = \det \widetilde{C} = m^{m-2}$.
    \end{proof}
	
Now we \bluer{state} a remarkable formula for the minors of the response matrix. 
Denote by $C_{I}^{J}$ the submatrix of the matrix $C$ made up of rows from a set $I$ and columns from a set $J$ in the natural order. Denote by \bluevarnew{$C_{X, Z}^{Y, Z}$} the block matrix 
\bluenew{
$\begin{pmatrix}
  \begin{matrix}
  C_{X}^{Y}
  \end{matrix}
  & \rvline & C_{X}^{Z} \\
\hline
  C_{Z}^{Y} & \rvline &
  \begin{matrix}
  C_{Z}^{Z}
  \end{matrix}
\end{pmatrix}.$} For example, with this notation, we have \bluenew{$\det C^{\{1\},\{2\}}_{\{3\},\{2\}}=-\det C^{\{1,2\}}_{\{2,3\}}$} because the latter matrix is obtained from the former by the swap of the rows. \bluenew{We will often use the notation $X := I \setminus J$, $Y := J \setminus I$, $Z := I \cap J$, $W := \{1, \dots, m \} \setminus (I \cup J)$. In this notation, $C_{X, Z}^{Y, Z} = \pm C_{I}^{J}$, where the sign depends on the number of row and column swaps required to arrange the indices into the natural order.} We set $\mathrm{sgn}(\pi)=-1$ for odd permutations $\pi$ and $\mathrm{sgn}(\pi)=+1$ for even ones (including the empty permutation).

For distinct 
vertices $x_{11}, \dots, x_{1n_1}, x_{21}, \dots, x_{2n_2}, \dots, x_{k1}, \dots, x_{kn_k}$,
denote by $$w(x_{11}, \dots, x_{1n_1}|x_{21}, \dots, x_{2n_2}|\dots|x_{k1}, \dots, x_{kn_k})$$ the sum of weights of all the spanning forests with exactly $k$ components such that:
\begin{itemize}
    \item $x_{11}, \dots, x_{1n_1}$ belong to the first connected component,
    
    $\vdots$
    
    \item $x_{k1}, \dots, x_{kn_k}$ belong to the $k$-th connected component,
    \item the other vertices belong to any connected component.
\end{itemize} 
\bluenew{Here we assume $n_1, \dots, n_k \ge 1$ but allow $k =0$; in the latter case, the sum is set to be~$0$.} 
For instance, the denominators in parts (1) and (2) of Theorem~\ref{Kirchhoff} are $w(1|2|\dots|m)$, and the numerators are 
$w(12\dots m)$ and
$w(ij|1|\dots|i-1|i+1|\dots|j-1|j+1|\dots|m)$ if $i<j$.

 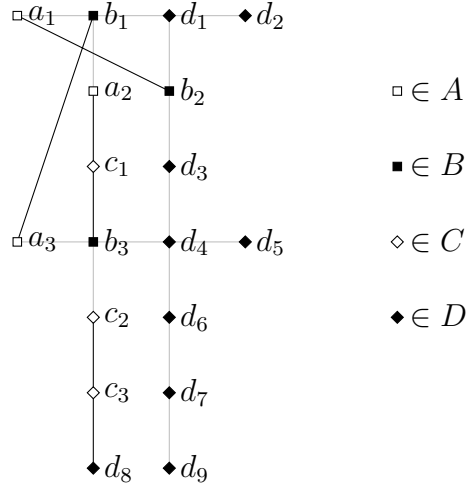
\begin{figure}[h]
\center
\begin{circuitikz} 
\draw [color = lightgray]     
               (1,2) to[short] (1,3)
               (1,6) to[short] (1,5)
               (0,6) to[short] (3,6)
               (2, 0) to (2,1 ) to (2,2) to (2, 3) to (2, 4) to (2, 5) to (2,6)
               (0,3) to (1,3) to (2,3) to (3,3);
             \draw
             (1,0) to (1,1) to (1,2)
             (0,3) to (1, 6)
             (1, 5) to (1, 4) to (1, 3)
             (0,6) to (2, 5)
             (1,0) node[diamondpole]{}
              (2,0) node[diamondpole]{}
              (1,1) node[odiamondpole]{}
              (2,1) node[diamondpole]{}
               (1,2) node[odiamondpole]{}
              (2,2) node[diamondpole]{}
               (1,3) node[squarepole]{}
              (2,3) node[diamondpole]{}
              (0,3) node[osquarepole]{}
               (1,4) node[odiamondpole]{}
              (2,4) node[diamondpole]{}
              (3,3) node[diamondpole]{}
               (1,5) node[osquarepole]{}
              (2,5) node[squarepole]{}
              (0,6) node[osquarepole]{}
              (1,6) node[squarepole]{}
              (2,6) node[diamondpole]{}
              (3,6) node[diamondpole]{};
            \draw {[anchor= west]
            
                            (1,0) node[]{$w_8$}
                            (2,0) node []{$w_9$}
                             (1,1) node[]{$z_3$}
                            (2,1) node []{$w_7$}
                            (1,2) node[]{$z_2$}
                            (2,2) node []{$w_6$}
                            (0,3) node[]{$x_3$}
                            (1,3) node[]{$y_3$}
                            (2,3) node []{$w_4$}
                            (3,3) node[]{$w_5$}
                            (1,4) node[]{$z_1$}
                            (2,4) node []{$w_3$}
                            (1,5) node[]{$x_2$}
                            (2,5) node []{$y_2$}
                            (0,6) node[]{$x_1$}
                            (1,6) node[]{$y_1$}
                            (2,6) node []{$w_1$}
                            (3,6) node[]{$w_2$}
                            };

                            \draw (5,2) node[diamondpole]{};
              \draw (5,3) node[odiamondpole]{};
              \draw (5,4) node[squarepole]{};
              \draw (5,5) node[osquarepole]{};
              \draw {[anchor= west]
              (5.02,5.02) node[]{ $\in X$}
              (5.02,4.02) node[]{ $\in Y$}
              (5.02,2.02) node[]{ $\in W$}
              (5.02,3.02) node[]{ $\in Z$}}; 
             \end{circuitikz}           
\caption{An electric network (shown in gray) with the boundary vertices decomposed into four sets $X,Y,Z,W$. Here $X=\{x_1,x_2,x_3\}$, $Y=\{y_1,y_2,y_3\}$, $Z=\{z_1,z_2,z_3\}$, \bluevarnew{and} $W=\{w_1,\dots,w_9\}$ are denoted by \bluevarnew{empty squares, filled squares, empty rhombi, and filled rhombi} respectively. 
The \bluevarnew{spanning} forest shown in black \bluevarnew{(including isolated vertices)} ``contributes'' to $C_{X, Z}^{Y, Z}$ with the sign $(-1)^{|X|}\mathrm{sgn}(\pi) =-1$. The vertices $x_1,x_2,x_3$ are joined with $y_2,y_3,y_1$ respectively by paths in the forest, hence the permutation
$\pi = \left( \begin{matrix}
			1 & 2 & 3\\
			2 & 3 & 1
		\end{matrix}\right)$,
 and the forest ``contributes'' to $w(x_1 y_2| x_2 y_3|x_3 y_1 |w_1|w_2|\dots|w_9)$. See Theorem~\ref{Kenyon-Wilson}.
} 
             \label{fig:KW}
             \end{figure}
             
\begin{theorem}[Kenyon--Wilson's all-minors matrix-tree theorem] \textup{(Cf.~\cite[Theorem~8.1]{KW09})} \label{Kenyon-Wilson}
    Decompose the set of boundary vertices of an electrical network into disjoint subsets \bluenew{$X, Y, Z, W$} such that $|X| = |Y|$. 
    Let $X = \{x_1, \dots, x_{|X|}\}$, $Y = \{y_1, \dots, y_{|Y|}\}$, $W=\{w_1, \dots, w_{|W|}\}$.
    Then 
	\begin{equation} \label{equation_Kenyon_Wilson} \det C_{X, Z}^{Y, Z} = \frac{(-1)^{|X|} \sum_{\pi} \mathrm{sgn}(\pi) \cdot w(
		x_1 y_{\pi(1)}| \dots |x_{|X|} y_{\pi(|X|)} | w_1| \dots | w_{|W|}
	)}{\sum_H w(H)},\end{equation} where the sum in the denominator is  over all valid forests $H$ in the electrical network, and the sum in the numerator is over all permutations $\pi \in S_{|X|}$ (see Fig.~\ref{fig:KW}).
	\end{theorem}

\bluenew{Here we do \emph{not} require that $x_1 < \dots < x_{|X|}$, $y_1 < \dots < y_{|Y|}$, $w_1 < \dots < w_{|W|}$.}

\begin{remark} Here we have corrected an overall sign error in the original statement \cite[Theorem 8.1]{KW09}. The sign error came from confusion between two common definitions of the response matrix, different by an overall sign. The proof of \cite[Theorem 8.1]{KW09} refers to \cite[Lemma A.1]{KW06}, where the sign convention for the response matrix was opposite (namely, \cite{KW06} uses the same definition as in our paper, while \cite{KW09} uses the opposite one). To fix this, the right side of~\cite[Eq.~(9)]{KW09} in the proof should be multiplied by $(-1) ^ {|Z|}$ because the minor has size $|Z|$. This results in the overall factor of $(-1) ^ {|Z|}$ in the final formula \cite[Theorem 8.1]{KW09}. To switch back to the notation of our paper, we need an additional factor of $(-1)^{|X|+|Z|}$ because minor~\eqref{equation_Kenyon_Wilson} has size $|X|+|Z|$. 
This results in a factor of $(-1)^{|Z|} (-1)^{|X|+|Z|} = (-1)^{|X|}$, which we added to \cite[Theorem 8.1]{KW09}. We have also checked the resulting formula in examples. \bluenew{For instance, for $X=\{i\}$, $Y=\{j\}$, and $Z$ empty, without the $(-1)^{|X|}$ factor, the result would contradict part (2) of Theorem~\ref{Kirchhoff}.}
\end{remark}

We conclude this section by recalling well-known formulae for the voltages and currents; cf.~\cite[Theorem~1.16]{Grimmett}. \emph{A spanning bitree} (or \emph{spanning thicket}) is a spanning forest with exactly two components.

\begin{proposition} \label{p-voltages-currents-electric-network}
Consider an electric circuit with two boundary vertices 
and prescribed voltages $U_1=1$ and $U_2=0$. 

(1) \bluenew{For each vertex $k$, we have}
$$U_k=\frac{\sum_{B_k} w(B_{k})}{\sum_{B} w(B)}=\frac{w(1k|2)}{w(1|2)},$$ where the sum in the denominator is over all spanning bitrees $B$ such that $1$ and $2$ are in different components, and the sum in the numerator is over all spanning bitrees $B_k$ such that $1$ and $2$ are in different components, whereas $k$ and $1$ are in the same component.

(2) \bluenew{For each oriented edge $kl$, we have}
$$I_{kl}=\frac{\sum_{T_{kl}} w(T_{kl})-\sum_{T_{lk}} w(T_{l k})}{\sum_{B} w(B)},$$ where the sum in the denominator is over all spanning bitrees $B$ such that $1$ and $2$ are in different components, and the sums in the numerator are over all spanning trees $T_{kl}$ and $T_{lk}$ such that the oriented path from $1$ to $2$ in the tree passes along the edge $kl$ 
in the direction from $k$ to $l$ and from $l$ to $k$ respectively.

\end{proposition}

\begin{proof}
It suffices to check that $I_{kl}$ and $U_k$ given by these expressions satisfy axioms (C) and (I). Let us verify axiom ~(C):
 \begin{align*}
 c_{kl}(U_k-U_l) 
 &=c_{kl}\frac{\sum w(B_{k}) -\sum w(B_{l}) }{\sum_{B} w(B)} 
= c_{kl}\frac{ w(1k|2) - w(1l|2) }{\sum_{B} w(B)}
\\
&= 
c_{kl}\frac{ w(1kl|2) +w(1k|2l) - w(1kl|2) - w(1l|2k)}{\sum_{B} w(B)} 
\\&= 
\frac{c_{kl} w(1k|2l) - c_{kl} w(1l|2k) }{\sum_{B} w(B)} 
\stackrel{(*)}{=}
\frac{\sum_{T_{kl}} w(T_{kl})-\sum_{T_{lk}} w(T_{lk})}{\sum_{B} w(B)}
=I_{kl}.
\end{align*}
Here equality~(*) holds because after adding the edge $kl$ to the spanning bitrees contributing to $w(1k|2l)$ and $w(1l|2k)$, the oriented path from $1$ to $2$ in the resulting spanning tree passes along the edge $kl$ 
in the direction from $k$ to $l$ and from $l$ to $k$ respectively, and 
the weight is multiplied by $c_{kl}$.
 
Let us verify axiom (I): \bluenew{for each vertex $k > 2$, we have}
$$\sum_{l=1}^{n} I_{kl}=\frac{\sum_{l=1}^{n}\left( \sum_{T_{kl}} w(T_{kl})-\sum_{T_{lk}} w(T_{l k}) \right)} {\sum_{B} w(B)}=0.$$
To prove the latter equality, let us show that each spanning tree contributes to the double sum exactly twice \bluenew{with opposite signs}. Take $2<k\le n$, $1\le i\le n$, and consider a tree $T_{ki}$ contributing to $\sum_{T_{kl}} w(T_{kl})$ for $l=i$. \bluenew{By the definition of the latter sum,} on the path from $1$ to $2$ in the tree $T_{ki}$, there are two edges with \bluenew{the} vertex $k$: $ki$ and $jk$ for some $j$, because $k>2$. Thus the same tree contributes to $\sum_{T_{lk}} w(T_{lk})$ for $l=j$ and the contributions cancel.

\bluenew{Finally, $ U_1 = 1$ and $U_2 = 0$ hold by construction and are obvious.}
\end{proof}

\section{Quick start: multiport networks}\label{sec-multiport}

Let us introduce our new concept in the simplest nontrivial particular case (see Fig.~\ref{fig:example}).

\begin{definition} \label{def-multiport}
    A \textit{$p$-port network} is a \bluenew{finite} connected graph with \bluenew{$p \ge 1$} marked disjoint ordered pairs of vertices (\emph{ports}) and a positive number (\textit{conductance}) assigned to each edge.

Fix enumeration of vertices $1,\dots,n$ such that the ports are $(1,2)$, \dots, $(m-1,m)$, where $m=2p$.
    
\textit{A $p$-port circuit} is a $p$-port network with $p$ real numbers (\textit{\bluenew{prescribed} voltage differences}) 
assigned to the $p$ ports. Denote by $\Delta U_{kl}$ the number assigned to a port $(k,l)$.

The \emph{voltages} $U_k$, where $1\le k\le n$, and the \emph{currents} $I_{kl}$, where $1\le k,l\le n$, in the $p$-port \bluenew{circuit} are determined (up to adding the same constant to all $U_k$) by axioms (C) and (I) above and the following two additional axioms:
 \begin{itemize}
		\item[(P)] \textit{Port isolation.} 
    	For each port $(k,l)$ we have $\sum_{i=1}^{n} I_{ki}+\sum_{i=1}^{n} I_{li}=0$.
		\item[(B)] \textit{Boundary conditions.} For each port $(k,l)$ we have $ U_k - U_l = \Delta{U_{kl}}$.
\end{itemize}
\end{definition}

\begin{theorem}\label{th-exist_uniq-multiport}
For any $p$-port circuit, there exists a unique collection of ${I_{kl}}$, ${1\le k,l\le n}$, and \bluenew{a unique (up to an additive constant)} collection of ${U_{k}}$, ${1\le k\le n}$, that satisfy conditions (C), (I), (P), (B).
\end{theorem}

This is a particular case of Theorem~\ref{exist_uniq} below, but let us present an alternative elementary proof due to O.~Angel (private communication).

\begin{proof}[Proof of Theorem~\ref{th-exist_uniq-multiport}]
    \emph{Existence}. First consider the case when $\Delta U_{12}=1$, but $\Delta U_{kl}=0$ for all the other ports~$(k,l)$. Set $U_1=1$ and $U_2=0$ to satisfy condition (B) for 
    port~$(1,2)$. Condition (B) for the other ports $(k,l)$ reads $U_k=U_l$. 
    Glue the two vertices $k$ and $l$ in each port $(k,l)\ne (1,2)$. Condition (P) for 
    port $(k,l)$ is equivalent to condition~(I) for the resulting combined vertex. Thus finding the voltages and currents in the $p$-port network reduces to finding voltages and currents in the electric network obtained by gluing and having only two boundary vertices $1$ and $2$. By Theorem~\ref{prelim_th}, \bluenew{namely, by the existence of a solution}, the case when $\Delta U_{12}=1$ and $\Delta U_{kl}=0$ for all ports $(k,l)\ne (1,2)$ follows.

    Analogously, the voltages and currents exist when some $\Delta U_{kl}=1$ and all the other \bluenew{prescribed voltage differences} 
    vanish. By linearity, they exist for arbitrary \bluenew{prescribed voltage differences}. 
    
    \emph{Uniqueness}. Suppose there are two collections of currents ${I^{I,II}_{kl}}$ and
voltages ${U^{I,II}_k}$ satisfying laws (C), (I), (P), (B) for some differences of prescribed voltages $\Delta U_{kl}^{I}$ = $\Delta U_{kl}^{II}$. Then the differences ${I_{kl}=I^{I}_{kl}-I^{II}_{kl}}$ and ${U_{k}=U^{I}_{k}-U^{II}_{k}}$ satisfy laws (C), (I), (P), (B) when $\Delta{U_{kl}}=0$ for all ports $(k,l)$. Gluing the two vertices $k$ and $l$ in each port $(k,l)\ne (1,2)$ we get an electric network with two boundary vertices $1$ and $2$. Since $U_1-U_2=\Delta U_{12}=0$, by Theorem~\ref{prelim_th},
\bluenew{namely, by the uniqueness of the solution},
it follows that all $U_{k}$ are the same and all $I_{kl}$ vanish. Hence
$U^{I}_{k}=U^{II}_{k}+\mathrm{const}$ and $I^{I}_{kl}=I^{II}_{kl}$ for all $k,l$.
\end{proof}

\textit{The response} of a $p$-port network is the linear map 
$$L\colon\mathbb{R}^{p}\to \mathbb{R}^{p}, \quad (\Delta{U_{12}},\Delta{U_{34}},\dots,\Delta{U_{2p-1,2p }}) \mapsto (I_{1}, I_{3},{\dots}, I_{2p-1})$$
taking the \bluenew{prescribed voltage differences} to the incoming currents.

In what follows we state our results in a more general setup but they are new and interesting already in the case of $p$-port networks, and one can restrict oneself to this case throughout.


\section{The definition of a superport network
}\label{sec-exist-unique}

    
    Let us introduce our new concept in full generality. Informally, in a superport network, boundary vertices are decomposed into several sets (superports). The voltages and currents are defined by the same laws (C) and (I) but the boundary conditions differ in two aspects. First, we impose an additional superport-isolation condition (P): the sum of the incoming currents is zero in each superport. Second, we relax the condition on voltages: instead of fixing all of them on the boundary, we fix voltage differences within each superport (this is referred to as a boundary condition (B)). 

    Let us summarize this concept in the form of a precise and self-contained definition.
        
    \begin{definition} \label{def-superport-network}
    \textit{A superport network with ${p}$ superports} is a \bluenew{finite} connected graph with \bluevarnew{$p \ge 1$} non-empty disjoint sets $A_1$, ${\dots}$, $A_p$ of vertices and a positive number (\textit{conductance}) assigned to each edge. The sets $A_1$, ${\dots}$, $A_p$ are called \textit{superports}. Their vertices are called \textit{boundary}, and the other vertices are called \textit{interior}. Denote by $M$ the set of boundary vertices and set $m := |M|$. For notational convenience, we assume that the graph has neither multiple edges nor loops.

    Fix an enumeration of the vertices ${1}$, ${2}$, ${\dots}$, ${n}$ such that ${1, \dots, m}$ are the boundary ones appearing in the order of the superports, i.e., $A_1 = \{1,2, \dots, |A_1|\}, A_2 = \{|A_1|+1, \dots, |A_1|+|A_2|\}$, and so on. Denote by ${kl}$ the edge between the vertices ${k}$ and ${l}$. Denote by $c_{kl}$ the conductance of the edge ${kl}$. Set $c_{kl} = 0$, if there is no edge between ${k}$ and ${l}$.

    The vertex with the maximal number in a superport is called a \textit{root}. The other boundary vertices are called \textit{non-root vertices}. Denote the set of all roots by $R = \{r_1, \dots, r_p\}$, where $r_k$ denotes the number of the root in the superport $A_k$, that is, $r_k := \sum \limits_{i=1}^k |A_i|$. Denote by $root(k)$ the root belonging to the same superport as a vertex $k$ (if $k$ is a root, then $root(k) = k$).

    \textit{A superport circuit} is a superport network with $m-p$ real numbers $\Delta{U_{kl}}$ (\textit{\bluenew{prescribed voltage differences}}) 
    assigned to $m-p$ non-root boundary vertices $k$, where $l = root(k)$.

Each superport circuit gives rise to certain numbers $U_k$, where $1\le k\le n$ (\textit{voltages} at the vertices), and $I_{kl}$ (\emph{currents} through the edges), where $1\le k,l\le n$.
    These numbers are defined by the following axioms:
    \begin{enumerate}
        \item[(C)] \textit{The Ohm law.} For each \bluenew{ordered pair of vertices $k, l$}, we have $I_{kl} = c_{kl}(U_k - U_l)$.
    	\item[(I)]\textit{The Kirchhoff current law.} For each vertex $k>m$ we have $\sum_{l=1}^{n} I_{kl}=0$.
        \item[(P)] \textit{Superport isolation.} 
        \item[(B)]\textit{Boundary conditions.} For any non-root boundary vertex $k$ and the root $l = root(k)$ we have $ U_k - U_l = \Delta{U_{kl}}$.
    \end{enumerate}
  The number $I_k := \sum_{l=1}^{n} I_{kl}$ is called \textit{the incoming current for vertex $k$}.
  \end{definition}

\begin{remark}
        For the last superport $A_p$, the equality ${\sum_{k\in A_p}{I_{k}=0}}$ also holds. Indeed, from axioms (C) and (I) it follows that the sum of all the incoming currents vanishes. Then from axiom (P) we obtain the equality ${\sum_{k\in A_i}{ I_{k}=0}}$ for all $i$ including $i=p$. We have excluded $i=p$ from axiom (P) to make the set of conditions minimal.
	\end{remark}
 
    \begin{example}
        If there is only one superport (i.e. $p = 1$), then we get an equivalent definition of an ordinary electrical circuit. If each superport has exactly $2$ vertices (i.e. $|A_k| = 2$ for $k=1,\dots,p$) then we get a $p$-port circuit.
    \end{example}
 
\begin{remark} \label{rem-superport-size-1}
  If some superport consists of just one vertex \bluenew{(and $p \neq 1$)}, then the superport can be just dropped and the vertex can be considered as an interior one. Indeed, this vertex is a root and does not contribute to boundary conditions (B), and the condition of superport isolation (P) is the same as Kirchhoff's law (I). Vice versa, any interior vertex can be replaced by a superport of size one.
	\end{remark}

\begin{remark}\label{another_motivation}
	 Let us provide yet another motivation for the notion of superport networks. Consider a network drawn in a planar domain without self-intersections so that the boundary vertices lie on the boundary of the domain. A \emph{conjugate} of the voltages is a function on the set of faces such that the difference of its values at two adjacent faces equals the current through their common edge (with an appropriate sign). If the domain is not simply-connected, e.g., an annulus, then this conjugate may not exist (or, if one prefers, may become multivalued). It is easy to see that the necessary and sufficient condition for the existence of a (single-valued) conjugate is vanishing of the total incoming current for each component of the boundary of the domain. This is exactly the superport condition. 
	\end{remark}	

\begin{theorem}\label{exist_uniq}
For any superport circuit, there exists a unique collection of ${I_{kl}}$, ${1\le k,l\le n}$, and a unique (up to an additive constant) collection of ${U_{k}}$, ${1\le k\le n}$, that satisfy conditions (C), (I), (P), (B).
\end{theorem}

The proof of Theorem~\ref{exist_uniq} is similar to \textup{\cite[Proof of Theorem~2.1]{PS}} and uses the following known lemma (alternatively, it can be proved analogously to Theorem~\ref{th-exist_uniq-multiport}).

\begin{lemma}[Energy conservation] \textup{~\cite[Lemma~5.1]{PS}} \label{energy} \label{fullpower} 
	Consider an electrical network with the vertices ${1, \dots, n}$ such that ${1, \dots, m}$ are the boundary ones.
    Assume that the numbers ${U_k}$, where ${1\le k\le n}$, and ${I_{kl}}$, where ${1\le k,l\le  n}$, satisfy conditions~\textup{(C)} and~\textup{(I)}. Then
	$$\sum\limits_{1\le k<l\le n}(U_k-U_l)I_{kl}
	= \sum\limits_{1\le u\le m}U_u I_u.$$
\end{lemma}

\begin{proof}[Proof of Theorem~\ref{exist_uniq}.] 	
	It is easy to see that adding the same constant to all the voltages ${U_{k}}$, where $1\le k\le n$, preserves all conditions (C), (I), (P), (B).
 We prove the existence and uniqueness of the collections under the additional condition $U_m=0$. 
 
\textit{Uniqueness}.
Suppose there are two collections of currents ${I^{I,II}_{kl}}$ and
voltages ${U^{I,II}_k}$ satisfying laws (C), (I), (P), (B), and ${U^{I}_m = U^{II}_m = 0}$ for some differences of prescribed voltages $\Delta U_{kl}^{I}$ = $\Delta U_{kl}^{II}$. Then the differences ${I_{kl}=I^{I}_{kl}-I^{II}_{kl}}$ and ${U_{k}=U^{I}_{k}-U^{II}_{k}}$ satisfy laws (C), (I), (P), (B) when $\Delta{U_{kl}}=0$ for all non-root vertices $k$, where $l = root(k)$. Thus, all voltages are the same within the same superport.

\bluenew{Adding all equations~(I), 
we get} ${\sum\limits_{k=1}^{m} I_k = 0.}$
Then by (P) we have ${\sum\limits_{k \in A_i}{ I_{k}=0}}$ for any $i$ including $i=p$. Then ${\sum\limits_{k \in A_i}{U_{k} I_{k}=0}}$ because $U_k = U_{root(k)}$ for all $k$. We get $\sum_{u=1}^m U_u I_u = 0$.

Then by (C) and Lemma~\ref{energy} we have $$\sum\limits_{1\le k<l\le n}(U_k-U_l)^2 c_{kl}
	=\sum\limits_{1\le k<l\le n}(U_k-U_l)I_{kl}
	= \sum\limits_{1\le u\le m}U_u I_u = 0.$$

For each ${k,l}$ we \bluenew{have ${c_{kl} \ge 0}$}. Thus each summand $ c_{kl}(U_k-U_l)^2=0$. Since the network is connected it follows that all the voltages ${U_{k}}$ are the same.
But $U_m = U_m^{I} - U_m^{II}=0$. Hence ${U_k=0}$, ${I_{kl}=0}$, and thus ${I^{I}_{kl}=I^{II}_{kl}}$, ${U^{I}_{k}=U^{II}_{k}}$ for each ${k,l}$.

\textit{Existence.} \bluenew{The number of equations in the system (C), (I), (P), and (B) equals the number of variables (if we fix $U_m=0$ and do not view $U_m$ as a variable). Indeed, condition (C) gives $n^2$ equations, (I) gives $n-m$ ones, (P) gives $p-1$ and (B) gives $m-p$. There are $n^2$ variables $I_{kl}$ and $n-1$ variables $U_k$ because $U_m=0$ is not viewed as a variable. Thus the number of equations and the number of variables are both $ n^2 + n - 1 $.}


We have proved that the system has a unique solution when $\Delta{U_{kl}}=0$ for all non-root boundary vertices $k$, where $l = root(k)$. By the finite-dimensional Fredholm alternative, it has a solution for any \bluenew{prescribed voltage differences} 
$\Delta{U_{kl}}$.
\end{proof}

\section{The response of a superport network} \label{sec-response}

\begin{definition}\label{def-response}
    \bluenew{Let $m>p$ and} let $x_1<\dots<x_{m-p}$ be the numbers \bluenew{assigned to} non-root vertices \bluevarnew{in the fixed enumeration}. The linear map 
    $$L\colon\mathbb{R}^{m-p}\to \mathbb{R}^{m-p}, \quad \left(\Delta{U_{x_1 root(x_1)}},\dots,\Delta{U_{x_{m-p} root(x_{m-p})}}\right) \mapsto (I_{x_1}, {\dots}, I_{x_{m-p}})$$ is \textit{the response} of the superport network. Denote by $L$ the matrix of the map as well.
\bluenew{The rows and columns of $L$ are enumerated by the set $M \setminus R=\{x_1, \dots, x_{m-p}\}$ of non-root vertices so that $L_i^j$ denotes the entry in the row $i \in M \setminus R$ and column $j \in M \setminus R$.}

    An electrical network is obtained from a superport network \emph{by unifying all the superports}, if both networks have the same graph, boundary vertices, and edge conductances. 
\end{definition}

Let us find out how the responses of the two networks are related. In particular, this will show that $L$ is symmetric.

\begin{theorem} \label{C-L} 
     Let $L$ and $C$ be the response matrices of a superport network  \bluenew{with $m > p$} and the electrical network obtained by unifying all the superports, respectively. Then $L$ is obtained from $C$ by Algorithm~\ref{alg:C2L}.
    %
\end{theorem}

\begin{algorithm}[tb]
    \caption{The C2L algorithm producing the response matrix of a superport network from the response matrix of the electrical network}
 \label{alg:C2L}
    \KwIn{A superport network and the response matrix $C$ of the electrical network obtained by unifying all the superports}
    \KwOut{The response matrix $L$ of the superport network}
    Remove the last row and the last column from $C$ and invert the resulting matrix\;
	For each non-root vertex $k$ with $root(k) \neq m$ subtract the column $root(k)$ from the column $k$\;
	For each non-root vertex $k$ with $root(k) \neq m$ subtract the row $root(k)$ from the row $k$\;
	For each root $k \neq m$ remove the row $k$ and the column $k$\;
    Invert the matrix and return the resulting matrix $L$.
    \end{algorithm}

\notdiplom{The idea of the proof is passing from the response, also known as the Dirichlet-to-Neuman map, to the Neuman-to-Dirichlet map, which is inverse in a sense. This makes the proof nearly obvious but still rather technical to write down.}

\begin{proof}[Proof of Theorem~\ref{C-L}]
    Step 1. Consider the composition of the response map of the electrical network restricted to the subspace $U_m=0$ and the projection to the subspace $I_m=0$. The matrix of the resulting map $(U_{1}, \dots, U_{m - 1}) \mapsto (I_{1}, \dots, I_{m - 1})$ is obtained from matrix $C$ by removing the last row and column. It follows from part (1) of Theorem~\ref{Kirchhoff} that the resulting matrix $\widetilde{C}$ is invertible. The inverse matrix defines the map $(I_{1}, \dots, I_{m-1}) \mapsto (U_{1}, \dots, U_{m-1})$, and it is the matrix obtained after step 1 of Algorithm~\ref{alg:C2L}.
	
	Step 2. Next, perform the following coordinate substitution in space $\mathbb{R}^{m-1}$: $$
	I'_k = 
	\begin{cases}
		I_k, &\text{if $k \neq root(k),$}\\
		\sum \limits_{l:root(l) = k} I_l, &\text{if $k = root(k).$}\\
	\end{cases}
	$$
	The matrix of the map $(I'_{1}, \dots, I'_{m-1}) \mapsto (U_{1}, \dots, U_{m-1})$ is obtained from the matrix of the map $(I_{1}, \dots, I_{m-1}) \mapsto (U_{1}, \dots, U_{m-1})$ by subtracting the column $root(k)$ from the column $k$ for all non-root vertices $k$ with $root(k)\neq m$, as in step 2 of Algorithm~\ref{alg:C2L}. Indeed, the columns of the former matrix are the image of the standard basis 
    $e_1',\dots,e_{m-1}'$ in $\mathbb{R}^{m-1}$, and our coordinate substitution takes $e_k'$ to itself if $k$ is a root and to $e_k'-e_{root(k)}'$ otherwise.
	
	Step~3. Now perform the coordinate substitution
 $$
	U'_k = 
	\begin{cases}
		U_k - U_{root(k)}, &\text{if $k \neq root(k)$},\\
		U_k, &\text{if $k = root(k)$},\\
	\end{cases}$$ 
 where we set $U_m:=0$.
 Then the map $(I'_{1}, \dots, I'_{m-1}) \mapsto (U'_{1}, \dots, U'_{m-1})$ is obtained from the map $ (I'_{1}, \dots, I'_{m-1}) \mapsto (U_{1}, \dots, U_{m-1})$ by subtracting the row $root(k)$ from the row $k$ for all non-root vertices $k$ with $root(k) \neq m$, as in step 3 of Algorithm~\ref{alg:C2L}.
 
Step 4. Consider the composition of the restriction of the map $(I'_{1}, \dots, I'_{m-1}) \mapsto (U'_{1}, \dots, U'_{m-1})$ to the subspace given by equations $I'_k=0$ for each root $k \neq m$, and the projection onto the subspace given by equations $U'_k = 0$ for each root $k \neq m$. The matrix of the resulting map $(I'_{x_1}, \dots, I'_{x_{m-p}}) \mapsto (U'_{x_1}, \dots, U'_{x_{m-p}})$ is obtained by removing row $k$ and column $k$ for each root $k \neq m$, as in step 4 of Algorithm~\ref{alg:C2L}. This map is $(I_{x_1}, {\dots}, I_{x_{m-p}}) \mapsto ({U_{x_1}-U_{root(x_1)}},\dots,{U_{x_{m-p}}-U_{root(x_{m-p})}})$, because $I'_{k} = I_{k}$ and $U'_{k} = U_k- U_{root(k)}$ for each non-root vertex $k$.

Step 5. Note that the composition of this map and the superport network response equals the identity map. Indeed, consider the map $(I'_{1}, \dots, I'_{m-1}) \mapsto (U_{1}, \dots, U_{m-1})$ introduced in step 2. Take any point $(U_{1}, \dots, U_{m-1})$ in the image of the restriction of this map to the subspace given by equations $I'_k=0$ for each root $k \neq m$.
Let the prescribed voltages of our electrical circuit be $(U_{1}, \dots, U_{m-1})$ and the differences of prescribed voltages in our superport circuit be $({U_{x_1}-U_{root(x_1)}},\dots,{U_{x_{m-p}}-U_{root(x_{m-p})}})$. Then the voltages and currents in the former circuit satisfy all conditions  (C),(I),(P),(B) in the latter circuit.
Indeed, axiom (B) holds since we have taken $\Delta{U_{x_i root(x_i)}} = U_{x_i} - U_{root(x_i)}$. Axioms (C) and (I) hold because they also hold in the electrical circuit. Axiom (P) holds because we consider the restriction to the subspace given by the equations $I'_k=0$ for each root $k \neq m$ and hence $\sum \limits_{l:root(l) = k} I_l = 0$. Thus, the composition in question is identical: it takes a vector of incoming currents $(I_{x_1}, {\dots}, I_{x_{m-p}})$ to itself. 

Since the matrix of the map $(I_{x_1}, {\dots}, I_{x_{m-p}}) \mapsto ({U_{x_1}-U_{root(x_1)}},\dots,{U_{x_{m-p}}-U_{root(x_{m-p})}})$, obtained after the 4 steps is square and $L$ is its left inverse, it follows that $L$ is the two-sided inverse, as in step 5 of Algorithm~\ref{alg:C2L}.
\end{proof}

Analogously, $L$ can be computed directly from the edge weights. The \emph{Kirchhoff matrix} of a superport network is the $n\times n$ matrix \bluenew{$K$} with the entries
$$
\bluevarnew{K_i^j}:=\begin{cases}
    -c_{ij}, &\text{if } i\ne j,\\
    \sum_{k=1}^n c_{ik}, &\text{if } i= j.
\end{cases}
$$ \bluenew{This definition is the same as for ordinary electrical networks.}


\begin{corollary} \label{cor-K2L}
    The response matrix $L$ of a superport network can be obtained from the Kirchhoff matrix $K$ using Algorithm~\ref{alg:C2L}, only at step 4 the rows and columns for all interior vertices are also removed.
\end{corollary}

\begin{proof}
    Declare each interior vertex to be a superport of size $1$ so that the vertex is a root. By Remark~\ref{rem-superport-size-1}, this does not affect the response matrix $L$. But the electric network obtained by unifying all the resulting superports has no interior vertices anymore, and its \bluevarnew{response} matrix is now exactly the Kirchhoff matrix $K$. Applying Theorem~\ref{C-L} for $C=K$, we get the desired corollary. Only step 4 of Algorithm~\ref{alg:C2L} is affected because the newly introduced superports have no non-root vertices.
\end{proof}

In particular, if the superport network has just one superport, then Algorithm~\ref{alg:C2L} is essentially the evaluation of $C$ in terms of $K$ using the Schur complement: \bluenew{inverting a matrix, then removing the same number of rows and columns, and inverting again results in the Schur complement of an appropriate block.} 

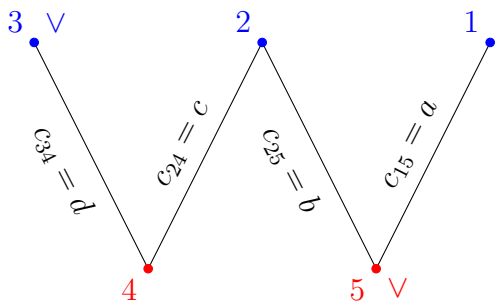
\begin{figure}[h!]
\center
\begin{circuitikz}
             \draw
             (0, 3) to[short, l_=${c_{34}=d}$] (1.5, 0)
             (1.5, 0) to[short, l=${c_{24}=c}$] (3, 3)
             (4.5, 0) to[short, l=${c_{25}=b}$] (3, 3)
             (4.5, 0) to[short, l=${c_{15}=a}$] (6, 3)
             (3,3) node[circ, color=blue]{}
             (0,3) node[circ, color=blue]{}
             (6,3) node[circ, color=blue]{}
             (1.5,0) node[circ, color=red]{}
             (4.5,0) node[circ, color=red]{}
            {[anchor=south east]
                            (0,3) node[color=blue] {3}
                            (6,3) node[color=blue] {1}
                            (3,3) node[color=blue] {2}
                            }
            {[anchor=north east]
                            (1.5,0) node[color=red] {4}
                           (4.5,0) node[color=red] {5}
                            }
            {[anchor=north west]
                           (4.5,0) node[color=red] {$\vee$}
                            }
              {[anchor=south west]
                            (0,3) node[color=blue] {$\vee$}
                            };
                            \end{circuitikz}
             \caption{A W-network. Colours depict superports and checkmarks indicate roots.}
             \label{fig:W} 
             \end{figure}

\begin{example}
A superport network with five vertices,
four edges $15, 25, 24, 34$ and two superports $\{1, 2, 3\}$ and $\{4, 5\}$ is called a \textit{W-network} (Fig.~\ref{fig:W}).
 Consider the W-network with the conductances
	$c_{15} = a, c_{25} = b, c_{24} = c, c_{34} = d$. The electrical network obtained by unifying all the superports has response matrix (equal to the Kirchhoff matrix)
	\begin{align*}\left( \begin{matrix}
			a & 0 & 0 & 0 & -a\\
			0 & b+c & 0 & -c & -b\\
			0 & 0 & d & -d & 0\\
			0 & -c & -d & c+d & 0\\
			-a & -b & 0 & 0 & a+b
		\end{matrix}\right).
	\end{align*}
 Let us calculate the response matrix of the superport W-network using Algorithm~\ref{alg:C2L}.
	\begin{enumerate}
		\item[Step 1.] Remove the last row and the last column from matrix $C$ and invert this matrix:
		\begin{align*}
			\left( \begin{matrix}
				a & 0 & 0 & 0\\
				0 & b+c & 0 & -c\\
				0 & 0 & d & -d\\
				0 & -c & -d & c+d
			\end{matrix}\right)^{-1}
			=
			\left( \begin{matrix}
				\frac{1}{a} & 0 & 0 & 0\\
				0 & \frac{1}{b} & \frac{1}{b} & \frac{1}{b}\\
				0 & \frac{1}{b} & \frac{bc+bd+cd}{bcd} & \frac{b+c}{bc}\\
				0 & \frac{1}{b} & \frac{b+c}{bc} & \frac{b+c}{bc}
			\end{matrix}\right).
		\end{align*}
		\item[Step 2.] Subtract column 3 from columns $1$ and $2$:
		\begin{align*}
			&\left( \begin{matrix}
				\frac{1}{a} & 0 & 0 & 0\\
				0 & \frac{1}{b} & \frac{1}{b} & \frac{1}{b}\\
				0 & \frac{1}{b} & \frac{bc+bd+cd}{bcd} & \frac{b+c}{bc}\\
				0 & \frac{1}{b} & \frac{b+c}{bc} & \frac{b+c}{bc}
			\end{matrix}\right)
			-
			\left( \begin{matrix}
				0 & 0 & 0 & 0\\
				\frac{1}{b} & \frac{1}{b} & 0 & 0\\
				\frac{bc+bd+cd}{bcd} & \frac{bc+bd+cd}{bcd} & 0 & 0\\
				\frac{b+c}{bc} & \frac{b+c}{bc} & 0 & 0
			\end{matrix}\right)
			= 
			\left( \begin{matrix}
				\frac{1}{a} & 0 & 0 & 0\\
				-\frac{1}{b} & 0 & \frac{1}{b} & \frac{1}{b}\\
				-\frac{bc+bd+cd}{bcd} & -\frac{c+d}{cd} & \frac{bc+bd+cd}{bcd} & \frac{b+c}{bc}\\
				-\frac{b+c}{bc} & -\frac{1}{c} & \frac{b+c}{bc} & \frac{b+c}{bc}
			\end{matrix}\right).
		\end{align*}
		\item[Step 3.] Subtract row $3$ from rows $1$ and $2$:
		\begin{align*}
			&\left( \begin{matrix}
				\frac{1}{a} & 0 & 0 & 0\\
				-\frac{1}{b} & 0 & \frac{1}{b} & \frac{1}{b}\\
				-\frac{bc+bd+cd}{bcd} & -\frac{c+d}{cd} & \frac{bc+bd+cd}{bcd} & \frac{b+c}{bc}\\
				-\frac{b+c}{bc} & -\frac{1}{c} & \frac{b+c}{bc} & \frac{b+c}{bc}
			\end{matrix}\right)
			-
			\left( \begin{matrix}
				-\frac{bc+bd+cd}{bcd} & -\frac{c+d}{cd} & \frac{bc+bd+cd}{bcd} & \frac{b+c}{bc}\\
				-\frac{bc+bd+cd}{bcd} & -\frac{c+d}{cd} & \frac{bc+bd+cd}{bcd} & \frac{b+c}{bc}\\
				0 & 0 & 0 & 0\\
				0 & 0 & 0 & 0
			\end{matrix}\right) = \\
			&\left( \begin{matrix}
				\frac{abc+abd+acd+bcd}{abcd} & \frac{c+d}{cd} & -\frac{bc+bd+cd}{bcd} & -\frac{b+c}{bc}\\
				\frac{c+d}{cd} & \frac{c+d}{cd} & -\frac{c+d}{cd} & -\frac{1}{c}\\
				-\frac{bc+bd+cd}{bcd} & -\frac{c+d}{cd} & \frac{bc+bd+cd}{bcd} & \frac{b+c}{bc}\\
				-\frac{b+c}{bc} & -\frac{1}{c} & \frac{b+c}{bc} & \frac{b+c}{bc}
			\end{matrix}\right).
		\end{align*}
		\item[Step 4.] Remove row $3$ and column $3$.
		\item[Step 5.] Invert the matrix and return the resulting matrix:
		\begin{align*}
			&L=\left( \begin{matrix}
				\frac{abc+abd+acd+bcd}{abcd} & \frac{c+d}{cd} & -\frac{b+c}{bc}\\
				\frac{c+d}{cd} & \frac{c+d}{cd} & -\frac{1}{c}\\
				-\frac{b+c}{bc} & -\frac{1}{c} & \frac{b+c}{bc}
			\end{matrix}\right)^{-1} \\&= \frac{1}{a+b+c+d}
			\left( \begin{matrix}
				ab+ac+ad & -ab-ac & ac+ad\\
				-ab-ac & ab+ac+bd+cd &  -ac+bd\\
				ac+ad & -ac+bd & ac+bc+ad+bd
			\end{matrix}\right).
		\end{align*}
	\end{enumerate}
 
 Notice that the entry $L^{2}_{4}=-ac+bd$ has terms of both signs (which is not possible for the entries of $C$). In particular, $\mathrm{sgn}(L^{2}_{4})$ can be arbitrary depending on $a,b,c,d$, and we observe no total positivity property as for planar electrical networks \bluenew{\cite[\S5.1]{CM}}.
\end{example}

\begin{remark}
    \bluenew{By definition, the entry $L^i_j$ can be interpreted as the \emph{incoming current at the vertex $i$ when the voltage at vertex $j$ is $1$, the voltages at the other vertices in the same superport are $0$, and the voltages within each of the remaining superports are equal}. Such an interpretation makes sense even if $i$ or $j$ are roots, and allows us to define an \emph{extended response matrix} of size $m\times m$. In particular, the entry $L^i_j$ does not depend on the choice of a root in each superport. 
    }
\end{remark}

\section{Voltages and currents in superport networks}
\label{sec-voltages-currents}

The matrix-tree theorem is another way to calculate the response of a network and its determinant. As a warm-up, we give combinatorial formulae for the response matrix, the voltages, and the currents in superport networks. Those formulae (but not the one for the determinant) are easily deduced from their analogs for ordinary electrical networks. 

To state those results, we generalize the notion of a valid forest to a superport network.

\begin{definition}\label{def-valid}
Introduce an equivalence relation on the set of vertices of a superport network: any two \bluenew{boundary} vertices are \emph{equivalent}, if and only if they are in the same superport, \bluenew{and} the equivalence class of each interior vertex consists of just a single vertex.
Denote by $\sim$ the resulting equivalence relation and by $[i]$ the equivalence class of a vertex~$i$.
A spanning forest in a superport network is \emph{valid}, if it becomes a spanning tree in the quotient by this equivalence relation (see Fig.~\ref{fig:fig-relatively-valid} to the bottom).

Notice that the validity of a forest depends on the choice of superports; in particular, a valid forest in a superport network may not be valid in the electrical network obtained by unifying all the superports (compare Fig.~\ref{fig:kirchhoff-matrix-tree} and \ref{fig:fig-relatively-valid}).

Now let $X$ be any subset of the set of vertices. Consider another equivalence relation: two \bluenew{boundary vertices not belonging to $X$} are \emph{$X$-equivalent} if and only if they are in the same superport, \bluenew{and the equivalence class of each interior vertex or a vertex of $X$ consists of just a single vertex.}
Denote still by $[i]$ the equivalence class of a vertex $i$.
A spanning forest is \emph{\bluevarnew{valid relative to} a vertex $i$} if it becomes a spanning tree in the quotient by the $\{i\}$-equivalence (see Fig.~\ref{fig:fig-relatively-valid} to the top).
The \emph{sign of a forest $G$}  (not necessarily valid) \emph{with respect to a pair of vertices $i$ and $j$} is 
$$\mathrm{sgn}(G,i,j)=
\begin{cases}
    +1, &\text{if either $i = j$, or $i$ and $j$ are in different components}\\
        &\text{of the quotient of $G$ by the $\{i,j\}$-equivalence},\\
    -1, &\text{otherwise.}
\end{cases}
$$
\end{definition}

\notdiplom{
\begin{theorem} \label{th-Lij} \textup{(See Fig.~\ref{fig:fig-relatively-valid})}
    The response matrix of a superport network has the entries
	$$L_{i}^j=\frac{\sum_G \mathrm{sgn}(G,i,j)\cdot w(G)}{\sum_F w(F)},$$
	where the sum in the denominator is over all valid forests $F$ and the sum in the numerator is over all forests $G$ \bluevarnew{valid relative to} each of the \bluenew{non-root} vertices $i$ and $j$ separately.
\end{theorem}

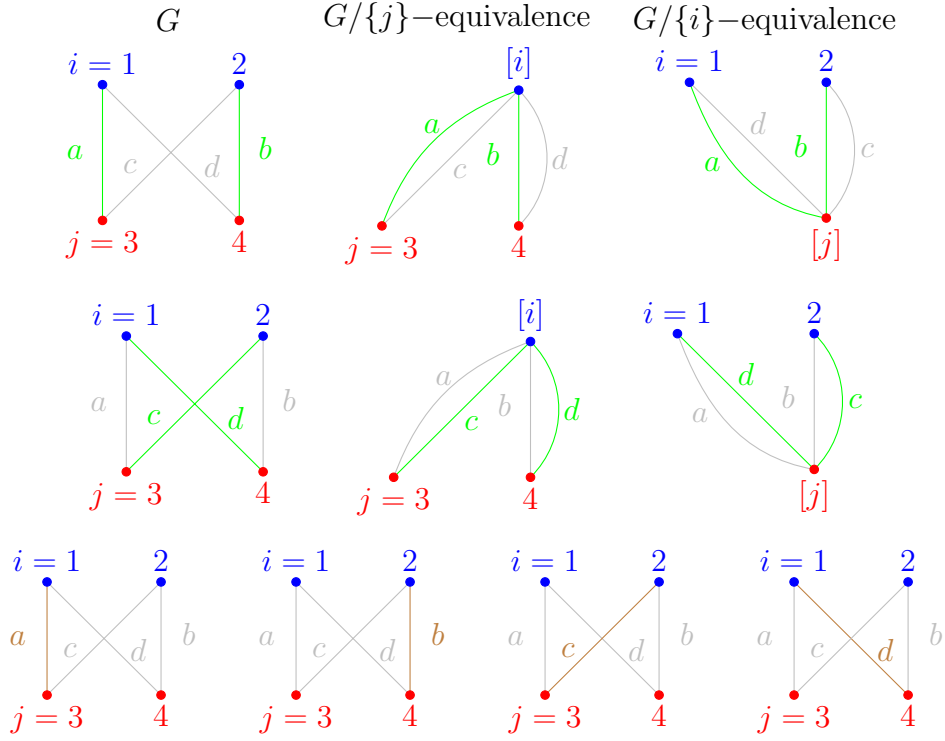
\begin{figure}[h]
    \centering

\begin{tabular}{c}
$G$\\

              \begin{circuitikz}[scale=0.6]
             \draw[color=green] (0,0) to[short, l=$a$] (0,3) 
             (3,3) to[short, l=$b$] (3, 0);
             \draw[color=lightgray]
             (0,3) to[] (3, 0)
             (0,0) to[] (3,3);
             \draw[color=lightgray] (0.6,1.2) node{$c$}
             (2.4,1.2) node{$d$};
             \draw {[anchor=north]
             (0,0) node[color=red]{$j=3$}};
             \draw {[anchor=south]
             (0,3) node[color=blue]{$i=1$}};
             \draw {[anchor=north]
             (3,0) node[color=red]{$4$}};
             \draw {[anchor=south]
             (3,3) node[color = blue]{$2$}};
             \draw (0,0)node[circ, color=red]{}
             (0,3)node[circ, color=blue]{}
             (3,0)node[circ, color=red]{}
             (3,3)node[circ, color=blue]{};
             \end{circuitikz}
\end{tabular} 
\begin{tabular}{c}
$G/\{j\}$-equivalence\\
              \begin{circuitikz}[scale=0.6]
             \draw [color=green] 
             (3,3) to[short, l_=$b$] (3, 0);
             \draw[color=lightgray] (3,3) to[out=-45, in=45] (3, 0);
             \draw[color=green] (0,0) to[out=70, in=-160] (3,3);
             \draw[color=lightgray] (0,0) to (3,3);
             \draw[color=lightgray] (1.7,1.2) node{$c$};
             \draw[color=green] (1.1,2.2) node{$a$};
             \draw[color=lightgray] (3.9,1.5) node{$d$};
             \draw {[anchor=north]
             (0,0) node[color=red]{$j=3$}};
             \draw {[anchor=north]
             (3,0) node[color=red]{$4$}};
             \draw {[anchor=south]
             (3,3) node[color = blue]{$[i]$}};
             \draw (0,0)node[circ, color=red]{}
             (3,0)node[circ, color=red]{};
             \draw (3,3)node[circ, color=blue]{};
             \end{circuitikz}
             \end{tabular} 
             \begin{tabular}{c}
             $G/\{i\}$-equivalence\\
             \begin{circuitikz}[scale=0.6]
             \draw[color=green]
             (3,3) to[short, l_=$b$] (3, 0);
             \draw[color=green]
             (0,3) to[out=-70, in=-190] (3,0);
             \draw[color=lightgray]
             (3,3) to[out=-45, in=45] (3, 0)
             (0,3) to[] (3, 0);
             \draw[color=green] (0.5,1.2) node{$a$};
             \draw[color=lightgray](1.5, 2.1) node{$d$};
             \draw[color=lightgray] (3.9,1.5) node{$c $};
             \draw {[anchor=south]
             (0,3) node[color=blue]{$i=1$}};
             \draw {[anchor=north]
             (3,0) node[color=red]{$[j]$}};
             \draw {[anchor=south]
             (3,3) node[color = blue]{$2$}};
             \draw 
             (0,3)node[circ, color=blue]{}
             (3,0)node[circ, color=red]{}
             (3,3)node[circ, color=blue]{};
             \end{circuitikz}
\end{tabular}

\begin{tabular}{c}
              \begin{circuitikz}[scale=0.6]
             \draw[color=lightgray] (0,0) to[short, l=$a$] (0,3) 
             (3,3) to[short, l=$b$] (3, 0);
             \draw[color=green]
             (0,3) to[] (3, 0)
             (0,0) to[] (3,3);
             \draw[color=green] (0.6,1.2) node{$c$}
             (2.4,1.2) node{$d$};
             \draw {[anchor=north]
             (0,0) node[color=red]{$j=3$}};
             \draw {[anchor=south]
             (0,3) node[color=blue]{$i=1$}};
             \draw {[anchor=north]
             (3,0) node[color=red]{$4$}};
             \draw {[anchor=south]
             (3,3) node[color = blue]{$2$}};
             \draw (0,0)node[circ, color=red]{}
             (0,3)node[circ, color=blue]{}
             (3,0)node[circ, color=red]{}
             (3,3)node[circ, color=blue]{};
             \end{circuitikz}
\end{tabular} 
\begin{tabular}{c}
              \begin{circuitikz}[scale=0.6]
             \draw [color=lightgray] 
             (3,3) to[short, l_=$b$] (3, 0);
             \draw[color=green] (3,3) to[out=-45, in=45] (3, 0);
             \draw[color=lightgray] (0,0) to[out=70, in=-160] (3,3);
             \draw[color=green] (0,0) to (3,3);
             \draw[color=green] (1.7,1.2) node{$c$};
             \draw[color=lightgray] (1.1,2.2) node{$a$};
             \draw[color=green] (3.9,1.5) node{$d$};
             \draw {[anchor=north]
             (0,0) node[color=red]{$j=3$}};
             \draw {[anchor=north]
             (3,0) node[color=red]{$4$}};
             \draw {[anchor=south]
             (3,3) node[color = blue]{$[i]$}};
             \draw (0,0)node[circ, color=red]{}
             (3,0)node[circ, color=red]{};
             \draw (3,3)node[circ, color=blue]{};
             \end{circuitikz}
             \end{tabular} 
             \begin{tabular}{c}
             \begin{circuitikz}[scale=0.6]
             \draw[color=lightgray]
             (3,3) to[short, l_=$b$] (3, 0);
             \draw[color=lightgray]
             (0,3) to[out=-70, in=-190] (3,0);
             \draw[color=green]
             (3,3) to[out=-45, in=45] (3, 0)
             (0,3) to[] (3, 0);
             \draw[color=lightgray] (0.5,1.2) node{$a$};
             \draw[color=green](1.5, 2.1) node{$d$};
             \draw[color=green] (3.9,1.5) node{$c $};
             \draw {[anchor=south]
             (0,3) node[color=blue]{$i=1$}};
             \draw {[anchor=north]
             (3,0) node[color=red]{$[j]$}};
             \draw {[anchor=south]
             (3,3) node[color = blue]{$2$}};
             \draw 
             (0,3)node[circ, color=blue]{}
             (3,0)node[circ, color=red]{}
             (3,3)node[circ, color=blue]{};
             \end{circuitikz}
\end{tabular}

\begin{tabular}{c}
\begin{tabular}{c}
              \begin{circuitikz}[scale=0.5]
             \draw[color=brown] (0,0) to[short, l=$a$] (0,3); 
             \draw[color=lightgray] (3,3) to[short, l=$b$] (3, 0);
             \draw[color=lightgray]
             (0,3) to[] (3, 0)
             (0,0) to[] (3,3);
             \draw[color=lightgray] (0.6,1.2) node{$c$}
             (2.4,1.2) node{$d$};
             \draw {[anchor=north]
             (0,0) node[color=red]{$j=3$}};
             \draw {[anchor=south]
             (0,3) node[color=blue]{$i=1$}};
             \draw {[anchor=north]
             (3,0) node[color=red]{$4$}};
             \draw {[anchor=south]
             (3,3) node[color = blue]{$2$}};
             \draw (0,0)node[circ, color=red]{}
             (0,3)node[circ, color=blue]{}
             (3,0)node[circ, color=red]{}
             (3,3)node[circ, color=blue]{};
             \draw[brown, thick] (3,3) circle (0.2);
             \draw[brown, thick] (3,0) circle (0.2);
             \end{circuitikz}
\end{tabular} 
\begin{tabular}{c}
              \begin{circuitikz}[scale=0.5]
             \draw[color=lightgray] (0,0) to[short, l=$a$] (0,3); 
             \draw[color=brown] (3,3) to[short, l=$b$] (3, 0);
             \draw[color=lightgray]
             (0,3) to[] (3, 0)
             (0,0) to[] (3,3);
             \draw[color=lightgray] (0.6,1.2) node{$c$}
             (2.4,1.2) node{$d$};
             \draw {[anchor=north]
             (0,0) node[color=red]{$j=3$}};
             \draw {[anchor=south]
             (0,3) node[color=blue]{$i=1$}};
             \draw {[anchor=north]
             (3,0) node[color=red]{$4$}};
             \draw {[anchor=south]
             (3,3) node[color = blue]{$2$}};
             \draw (0,0)node[circ, color=red]{}
             (0,3)node[circ, color=blue]{}
             (3,0)node[circ, color=red]{}
             (3,3)node[circ, color=blue]{};
             \draw[brown, thick] (0,3) circle (0.2);
             \draw[brown, thick] (0,0) circle (0.2);
             \end{circuitikz}
\end{tabular}
\begin{tabular}{c}
              \begin{circuitikz}[scale=0.5]
             \draw[color=lightgray] (0,0) to[short, l=$a$] (0,3); 
             \draw[color=lightgray] (3,3) to[short, l=$b$] (3, 0);
             \draw[color=lightgray]
             (0,3) to[] (3, 0);
             \draw[color=brown] (0,0) to[] (3,3);
             \draw[color=brown] (0.6,1.2) node{$c$};
             \draw[color=lightgray] (2.4,1.2) node{$d$};
             \draw {[anchor=north]
             (0,0) node[color=red]{$j=3$}};
             \draw {[anchor=south]
             (0,3) node[color=blue]{$i=1$}};
             \draw {[anchor=north]
             (3,0) node[color=red]{$4$}};
             \draw {[anchor=south]
             (3,3) node[color = blue]{$2$}};
             \draw (0,0)node[circ, color=red]{}
             (0,3)node[circ, color=blue]{}
             (3,0)node[circ, color=red]{}
             (3,3)node[circ, color=blue]{};
             \draw[brown, thick] (0,3) circle (0.2);
             \draw[brown, thick] (3,0) circle (0.2);
             \end{circuitikz}
\end{tabular} 
\begin{tabular}{c}
              \begin{circuitikz}[scale=0.5]
             \draw[color=lightgray] (0,0) to[short, l=$a$] (0,3); 
             \draw[color=lightgray] (3,3) to[short, l=$b$] (3, 0);
             \draw[color=brown]
             (0,3) to[] (3, 0);
             \draw[color=lightgray] (0,0) to[] (3,3);
             \draw[color=lightgray] (0.6,1.2) node{$c$};
             \draw[color=brown] (2.4,1.2) node{$d$};
             \draw {[anchor=north]
             (0,0) node[color=red]{$j=3$}};
             \draw {[anchor=south]
             (0,3) node[color=blue]{$i=1$}};
             \draw {[anchor=north]
             (3,0) node[color=red]{$4$}};
             \draw {[anchor=south]
             (3,3) node[color = blue]{$2$}};
             \draw (0,0)node[circ, color=red]{}
             (0,3)node[circ, color=blue]{}
             (3,0)node[circ, color=red]{}
             (3,3)node[circ, color=blue]{};
             \draw[brown, thick] (3,3) circle (0.2);
             \draw[brown, thick] (0,0) circle (0.2);
             \end{circuitikz}
\end{tabular} 
\end{tabular}
$$L_1^3 = \frac{\begingroup\color{green}cd-ab\endgroup}
    {\begingroup\color{brown}a+b+c+d\endgroup}$$\\
\caption{The matrix-tree theorem for superport networks. Superports are shown in red and blue. Forests $G$ \bluevarnew{valid relative to} each of the vertices $i$ and $j$ separately (and their quotients by the $\{i\}$- and $\{j\}$-equivalences) are shown in green \bluevarnew{(see the first two rows)}. Valid forests are shown in brown \bluevarnew{(see the last row)}. The rest of the edges are shown in light gray. See Definition~\ref{def-valid} and Theorem~\ref{th-Lij}. \bluevarnew{In this example, the $\{i,j\}$-equivalence is trivial, that is, the quotient of $G$ by the $\{i,j\}$-equivalence equals $G$.}}
\label{fig:fig-relatively-valid} 
\end{figure}

The idea of the proof is as follows. By definition, $L_{i}^j$ is the current $I_j$ when $\Delta U_{i, root(i)} = 1$ and $\Delta U_{k,root(k)}=0$ for all $k\ne i$. By axiom~(B), the latter implies that all the vertices $k\ne i$ within each port have the same voltages and can be considered equivalent. This is exactly the $\{i\}$-equivalence. The quotient of the superport network 
by the $\{i\}$-equivalence can be viewed as an electrical network with two boundary vertices $[i]$ and $[root(i)]$: each superport not containing $i$ is glued into a single vertex and thus dropped; see Remark~\ref{rem-superport-size-1}. 
 
\begin{lemma} \label{gluing}
Let a superport circuit have prescribed voltage differences $\Delta U_{i, root(i)} = 1$ for some non-root vertex $i$ and $\Delta U_{k,root(k)}=0$ for all $k\ne i$. Let the quotient of the superport circuit by the $\{i\}$-equivalence be viewed as an electrical circuit with two boundary vertices $[i]$ and $[root(i)]$ and prescribed voltages 
$U_{[i]}=1$ and $U_{[root(i)]}=0$. 
Then in both circuits, the currents are the same and the voltages are the same up to adding a constant.
\end{lemma}

\begin{proof} Consider the currents and voltages in the superport circuit.
They satisfy axioms (C), (I), (P), (B) and
give rise to well-defined currents and voltages in the quotient electrical circuit because we glue together vertices with the same voltage only.
Let us check that the resulting currents and voltages satisfy axioms (C) and (I). Axiom~(C) remains literally the same. Axiom~(I) needs to be checked only for the newly appeared interior vertices in the quotient. Each such vertex is obtained from a superport not containing $i$. The sum of currents coming out of the vertex equals the sum of currents coming out of the corresponding superport, vanishing by axiom (P). 
\end{proof}

\begin{remark} \label{gluing_remark}
An analogous assertion holds for the quotient by the $\{i, j\}$-equivalence.
If $i\ne j$ belong to the same superport, then the resulting electrical circuit has three boundary vertices $[i]$, $[j]$, and $[root(i)]$ and the prescribed voltages $U_{[i]}=1$, $U_{[j]}=U_{[root(i)]}=0$. 
\end{remark}


The proof of Theorem~\ref{th-Lij} is a simple consequence of \bluevarnew{Lemma~\ref{gluing}} and Proposition~\ref{p-voltages-currents-electric-network}.

\begin{proof}[Proof of Theorem~\ref{th-Lij}]
\bluenew{Take two non-root vertices $i$ and $j$}. Set $\Delta U_{i, root(i)} = 1$ and $\Delta U_{k,root(k)}=0$ for all $k\ne i$ so that $L_{i}^j=I_j$. Consider the following two cases.

\emph{Case 1}: $i$ and $j$ are in the same superport. 
Take the quotient by the $\{i, j\}$-equivalence \bluenew{and denote by $[x]$ the equivalence class of a vertex $x$ (throughout Case 1).} By Remark~\ref{gluing_remark}, the currents and the voltages in the resulting electric circuit remain the same. 
For $i\ne j$, express $L_{i}^j=I_j=C_{[i]}^{[j]}$ using part (2) of Theorem~\ref{Kirchhoff}. 
The resulting sums over forests $H$ and $G$ are exactly the ones in \bluenew{Theorem~\ref{th-Lij}} after taking the quotient. Indeed, take any forest $G$ as in \bluenew{Theorem~\ref{th-Lij}} and consider its quotient \bluenew{$G'$. 
Note that the quotient of $G$ by the $\{i\}$-equivalence is obtained from the quotient $G'$ by the $\{i, j\}$-equivalence by gluing the component containing $[j]$ with the component containing $[root(i)]$.} Since $G$ is \bluevarnew{valid relative to} $i$, it follows that either
$[i]$ and $[j]$ or $[i]$ and $[root(i)]$ lie in the same component \bluenew{of $G'$}. Since $G$ is \bluevarnew{valid relative to} $j$, it follows that either
$[i]$ and $[j]$ or $[j]$ and $[root(j)]=[root(i)]$ lie in the same component \bluenew{of $G'$}.   
Hence $[i]$ and $[j]$ lie in the same component \bluenew{of $G'$} and $\mathrm{sgn}(G,i,j)=-1$. 

For $i=j$, we express $L_{j}^j = C_{[j]}^{[j]} = -C_{[j]}^{[root(j)]}$ analogously; the latter equality holds because there are only two boundary vertices $[j]=[i]$ and $[root(j)]$.


\emph{Case 2}: $i$ and $j$ are in different superports. \bluenew{This time} take the quotient by the $\{i\}$-equivalence \bluenew{and denote by $[x]$ the equivalence class of a vertex $x$ (throughout Case~2).} By Lemma~\ref{gluing}, the currents and the voltages in the resulting electric circuit remain the same. We use the same notation for the edges before and after taking the quotient.

By Proposition~\ref{p-voltages-currents-electric-network} we get
\begin{align*} L_{i}^j=I_{j} = \sum_{k=1}^{n} I_{jk}&=\frac{\sum_{k=1}^{n}\left( \sum_{F_{jk}} w(F_{jk})-\sum_{F_{kj}} w(F_{kj}) \right)}{\sum_{F} w(F)},
\end{align*} 
where the sums are over forests $F$, $F_{jk}$, and $F_{k j}$ in the superport network such that their quotients $B$, $T_{jk}$, and $T_{k j}$ \bluenew{respectively} 
satisfy the \bluevarnew{conditions} listed in Proposition~\ref{p-voltages-currents-electric-network} \bluenew{with vertices $1$ and $2$ replaced with $[i]$ and $[root(i)]$.}



The denominator in the right side is the same as in the theorem: 
\bluenew{$[i]$ and $[root(i)]$} lie in distinct components of the bitree $B$ if and only if $B/\sim$ 
is a spanning tree, that is, $F$ is valid.

Consider a forest $F_{jk}$ contributing to the numerator. The quotient  $T_{jk}$ 
contains a unique oriented path from vertex $[i]$ to $[root(i)]$ which passes along the edge $jk$ and also another edge $uv$ with $[v]=[j]$, because $[j]\ne [i]$. 

If $v=j$ then the same forest $F_{jk}$ contributes also to $\sum_{F_{kj}} w(F_{kj})$ for $k=u$ and the contributions cancel.


If $v\ne j$ 
then the forest $F_{jk}$ contributes to the numerator exactly once with positive sign. \bluenew{Let us show that $F_{jk}$ is valid relative to $i$ and $j$ separately and $\mathrm{sgn}(F_{jk},i,j)=+1$.} Since the quotient $T_{jk}$ \bluenew{of $F_{jk}$} by the $\{i\}$-equivalence is a spanning tree, it follows that $F_{jk}$ is \bluevarnew{valid relative to} $i$. The quotient by $\{j\}$-equivalence is still a spanning tree: gluing the vertices $[i]$ and $[root(i)]$ together does not add a cycle because the path between them is split at the point $[j]$. Thus $F_{jk}$ is \bluevarnew{valid relative to} $j$ as well. The quotient \bluenew{of $F_{jk}$} by the $\{i,j\}$-equivalence is a spanning bitree with $i$ and $j$ in different components ($uv$ and $i$ lie in one component, whereas $jk$ and $root(i)$ lie in the other one). Thus $\mathrm{sgn}(F_{jk},i,j)=+1$. Conversely, each forest $G$ 
in the theorem with $\mathrm{sgn}(G,i,j)=+1$ equals $F_{jk}$ for some~$k$.

Similarly, the sum over forests $F_{kj}$ in the numerator reduces to the sum over
forests $G$ in the theorem with $\mathrm{sgn}(G,i,j)=-1$. We arrive at the desired expression.
\end{proof}
}

Using this approach, we can find the currents and the voltages in a superport circuit.

\begin{proposition} \label{p-voltages-currents-superport-network}
\bluenew{Let} a superport circuit 
\bluevarnew{have} prescribed voltage differences $\Delta U_{ij}$.
Then

(1) \bluenew{for each vertex $k$, we have}
$$U_k= \sum_{i \in M \setminus R} \Delta U_{i, root(i)}  \frac{\sum_{F_k} w(F_k)}{\sum_{F} w(F)}+\mathrm{const},$$ where the sum in the denominator is over all valid forests $F$, and the sum in the numerator is over all valid forests $F_k$ such that \bluenew{$[k]$} and \bluevarnew{$[i]$} are in the same component in the quotient of the forest by the $\{i\}$-equivalence.

(2) \bluenew{for each oriented edge $kl$, we have}
$$I_{kl}= \sum_{i \in M \setminus R} \Delta U_{i, root(i)} \frac{\sum_{F_{kl}} w(F_{kl})-\sum_{F_{l k}} w(F_{l k })}{\sum_{F} w(F)},$$ where the sum in the denominator is over all valid forests $F$, and the sums in the numerator are over all forests $F_{kl}$ and $F_{lk}$ \bluenew{valid relative to} the vertex $i$ such that the oriented path from \bluevarnew{$[i]$} to \bluevarnew{$[root(i)]$} in the quotient of the forest by the $\{i\}$-equivalence passes along the edge $kl$ in the direction from $k$ to $l$ and from $l$ to $k$ respectively.

\end{proposition}

\begin{proof}
By the linearity, it suffices to prove the proposition in the case when $\Delta U_{i, root(i)} = 1$ for some non-root vertex $i$ and $\Delta U_{k,root(k)}=0$ for all $k\ne i$. 
Take the quotient by the $\{i\}$-equivalence. By Lemma~\ref{gluing}, the currents and voltages in the resulting electrical circuit remain the same. 
They are given by expressions from Proposition~\ref{p-voltages-currents-electric-network}. The resulting sums over trees and bitrees are exactly the ones in the proposition after taking the quotient. 
\end{proof}


\section{Matrix-tree theorem for superport networks}
\label{sec-matrix-tree}

Now we state our main result: the determinantal matrix-tree theorem for superport networks. The statement relies on Definitions~\ref{def-superport-network}, \ref{def-response}, and \ref{def-valid} only. Although the statement is similar to the one for electrical networks (Theorem~\ref{Kirchhoff}), the proof is much more involved and contains new ideas.

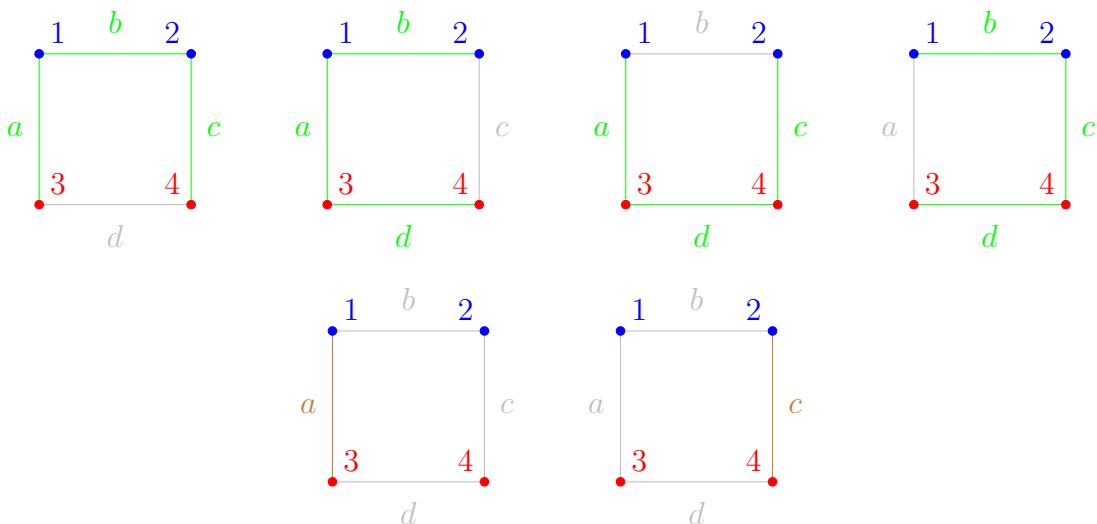
\begin{figure}[h]
    \centering

\begin{tabular}{cl}  
         \begin{tabular}{c}
           \begin{circuitikz}
           \draw [color = green] (0,0) to[short, l=$a$] (0,2)
                (0,2) to[short, l=$b$] (2,2)
               (2,2) to[short, l=$c$] (2,0);
                 \draw [color = lightgray]     
               (2,0) to[short, l=$d$] (0,0)
                 (0,0) node[circ, color=red]{}
                 (0,2) node[circ, color=blue]{}
                 (2,0) node[circ, color=red]{}
                 (2,2) node[circ, color=blue]{}
                
                {[anchor=south east]
                (2,0) node[color=red] {4}
                (2,2) node[color=blue] {2}
                }
                {[anchor=south west]
                (0,0) node[color=red] {3}
                (0,2) node[color=blue] {1}
                };
            \end{circuitikz}
         \end{tabular}
           & \begin{tabular}{l}

             \begin{circuitikz}
           \draw [color = green] (0,0) to[short, l=$a$] (0,2)
                (0,2) to[short, l=$b$] (2,2)
                (2,0) to[short, l=$d$] (0,0);
                 \draw [color = lightgray]
               (2,2) to[short, l=$c$] (2,0)
                 (0,0) node[circ, color=red]{}
                 (0,2) node[circ, color=blue]{}
                 (2,0) node[circ, color=red]{}
                 (2,2) node[circ, color=blue]{}
                
                {[anchor=south east]
                (2,0) node[color=red] {4}
                (2,2) node[color=blue] {2}
                }
                {[anchor=south west]
                (0,0) node[color=red] {3}
                (0,2) node[color=blue] {1}
                };
            \end{circuitikz}
         \end{tabular}  
\end{tabular}
\begin{tabular}{cl}  
         \begin{tabular}{c}
           \begin{circuitikz}
           \draw [color = green] (0,0) to[short, l=$a$] (0,2)
                (2,0) to[short, l=$d$] (0,0)
               (2,2) to[short, l=$c$] (2,0);
                 \draw [color = lightgray]
               (0,2) to[short, l=$b$] (2,2)
                 (0,0) node[circ, color=red]{}
                 (0,2) node[circ, color=blue]{}
                 (2,0) node[circ, color=red]{}
                 (2,2) node[circ, color=blue]{}
                
                {[anchor=south east]
                (2,0) node[color=red] {4}
                (2,2) node[color=blue] {2}
                }
                {[anchor=south west]
                (0,0) node[color=red] {3}
                (0,2) node[color=blue] {1}
                };
            \end{circuitikz}
         \end{tabular}
           & \begin{tabular}{l}

             \begin{circuitikz}
           \draw [color = green] 
                (0,2) to[short, l=$b$] (2,2)
                (2,0) to[short, l=$d$] (0,0)
               (2,2) to[short, l=$c$] (2,0);
                 \draw [color = lightgray]
                 (0,0) to[short, l=$a$] (0,2)
                 (0,0) node[circ, color=red]{}
                 (0,2) node[circ, color=blue]{}
                 (2,0) node[circ, color=red]{}
                 (2,2) node[circ, color=blue]{}
                
                {[anchor=south east]
                (2,0) node[color=red] {4}
                (2,2) node[color=blue] {2}
                }
                {[anchor=south west]
                (0,0) node[color=red] {3}
                (0,2) node[color=blue] {1}
                };
            \end{circuitikz}
         \end{tabular} 
\end{tabular}   

\begin{tabular}{cl}  
         \begin{tabular}{c}
           \begin{circuitikz}
           \draw [color = brown] (0,0) to[short, l=$a$] (0,2);
                 \draw [color = lightgray]
                 (2,0) to[short, l=$d$] (0,0)
               (2,2) to[short, l=$c$] (2,0)
               (0,2) to[short, l=$b$] (2,2)
                 (0,0) node[circ, color=red]{}
                 (0,2) node[circ, color=blue]{}
                 (2,0) node[circ, color=red]{}
                 (2,2) node[circ, color=blue]{}
                
                {[anchor=south east]
                (2,0) node[color=red] {4}
                (2,2) node[color=blue] {2}
                }
                {[anchor=south west]
                (0,0) node[color=red] {3}
                (0,2) node[color=blue] {1}
                };
                \draw[brown, thick] (2,2) circle (0.15);
                \draw[brown, thick] (2,0) circle (0.15);
            \end{circuitikz}
         \end{tabular}
           & \begin{tabular}{l}

             \begin{circuitikz}
           \draw [color = brown] 
               (2,2) to[short, l=$c$] (2,0);
                 \draw [color = lightgray]
                 (0,2) to[short, l=$b$] (2,2)
                (2,0) to[short, l=$d$] (0,0)
                 (0,0) to[short, l=$a$] (0,2)
                 (0,0) node[circ, color=red]{}
                 (0,2) node[circ, color=blue]{}
                 (2,0) node[circ, color=red]{}
                 (2,2) node[circ, color=blue]{}
                
                {[anchor=south east]
                (2,0) node[color=red] {4}
                (2,2) node[color=blue] {2}
                }
                {[anchor=south west]
                (0,0) node[color=red] {3}
                (0,2) node[color=blue] {1}
                };
                \draw[brown, thick] (0,2) circle (0.15);
                \draw[brown, thick] (0,0) circle (0.15);
            \end{circuitikz}  
         \end{tabular}  \\
\end{tabular}
$$L = \left( \begin{matrix}
    \frac{ab+bc+ac}{a+c} & \frac{-ac}{a+c} \\
    \frac{-ac}{a+c} & \frac{ac + ad+cd}{a+c}
\end{matrix}\right) \quad \quad \quad  \det L =\frac{\begingroup\color{green}abc+abd+acd+bcd\endgroup}
{\begingroup\color{brown}a+c\endgroup}$$

    \caption{The determinantal matrix-tree theorem for superport networks. Superports are shown in red and blue. Spanning trees and valid forests are shown in green and brown, respectively. The rest of the edges are shown in light gray. See Theorem~\ref{main}.}  
    \label{fig:superport-matrix-tree}
\end{figure}

\begin{theorem}[Determinantal matrix-tree theorem for superport networks] \label{main}
For any superport network \bluenew{with $m > p$} the determinant of the response matrix $L$ equals
	$$\det L = \frac{\sum_T w(T)}{\sum_F w(F)},$$
where the sum in the denominator is over all valid forests $F$ in the superport network, and the sum in the numerator is over all spanning trees $T$ (see Fig.~\ref{fig:superport-matrix-tree}).
\end{theorem}

\bluenew{Let us illustrate this result by a known combinatorial corollary 
\cite[Problem~2.2.3(3)]{Chernov-etl-16}.}

\begin{corollary}[Generalized Cayley's formula] \label{generalized_Cayley}
Let $T_1, \dots, T_r$ be trees with the sets of vertices $A_1, \dots, A_r$ respectively. The number of trees whose set of vertices is the disjoint union $A_1\sqcup \dots \sqcup A_r$ and which contain $T_1 \sqcup \dots \sqcup T_r$ is equal to 
$\left( \sum \limits_{i=1}^r |A_i| \right) ^ {r-2} \cdot \prod \limits_{i=1}^r |A_i|.$
\end{corollary}

\begin{proof}
Let $n := \sum \limits_{i=1}^r |A_i|.$
Consider the superport network on the complete graph on $n$ vertices with unit edge conductances and with $r$ superports $A_1, \dots, A_r$.

The response matrix for the electrical network obtained by unifying all the superports is given by
$C_{i}^{i} = n-1$ for any $i$ and $C_{i}^{j}=-1$ for any $i \neq j$.

Convert this $C$ to the response matrix $L$ of the superport network using Algorithm~\ref{alg:C2L}.

After step 1, we obtain the matrix $\widetilde{C}^{-1}$. It is easily checked that $(\widetilde{C}^{-1})_{i}^{i} = 2/n$, $(\widetilde{C}^{-1})_{i}^{j}=1/n$ for $i \neq j$ by multiplying $\widetilde{C}$ and $\widetilde{C}^{-1}$.

After steps \bluenew{2--4}, we get a block-diagonal matrix with block sizes $|A_1| - 1, \dots, |A_r|-1$, elements on the diagonal equal to $2/n$, and elements not on the diagonal equal to $1/n$ in each block. \bluer{For a block of size $|A_i|-1$, the determinant is calculated as follows:

\begin{align*}
\det \left(\begin{array}{ccccc}
    \frac{2}{n} & \frac{1}{n} & \ldots & \frac{1}{n} \\
    \frac{1}{n} & \frac{2}{n} & \ldots & \frac{1}{n} \\
    \vdots & \vdots & \ddots & \vdots \\
    \frac{1}{n} & \frac{1}{n} & \ldots & \frac{2}{n}
\end{array}\right) 
&= \frac{1}{n^{|A_i|-1}} \det \left(\begin{array}{cccc}
    2 & 1 & \ldots & 1 \\
    1 & 2 & \ldots & 1 \\
    \vdots & \vdots & \ddots & \vdots \\
    1 & 1 & \ldots & 2
\end{array}\right) 
= \frac{1}{n^{|A_i|-1}} \det \left(\begin{array}{cccc}
    |A_i| & |A_i| & \ldots & |A_i| \\
    1 & 2 & \ldots & 1 \\
    \vdots & \vdots & \ddots & \vdots \\
    1 & 1 & \ldots & 2
\end{array}\right) \\
&= \frac{|A_i|}{n^{|A_i|-1}} \det\left(\begin{array}{cccc}
    1 & 1 & \ldots & 1 \\
    1 & 2 & \ldots & 1 \\
    \vdots & \vdots & \ddots & \vdots \\
    1 & 1 & \ldots & 2
\end{array}\right) = \frac{|A_i|}{n^{|A_i|-1}} \det\left(\begin{array}{cccc}
    1 & 1 & \ldots & 1 \\
    0 & 1 & \ldots & 0 \\
    \vdots & \vdots & \ddots & \vdots \\
    0 & 0 & \ldots & 1
\end{array}\right) = \frac{|A_i|}{n^{|A_i|-1}}.
\end{align*}}
\bluenew{Another method to find the determinant of this matrix involves examining its spectrum. The spectrum of the matrix with all elements equal to $1$ is $\{|A_i| - 1, 0, 0, \ldots, 0\}$. Adding the identity matrix shifts the eigenvalues to $\{|A_i|, 1, 1, \ldots, 1\}$. Thus, we can compute the determinant of the block by factoring out ${1}/{n^{|A_i|-1}}$, yielding the same result.}

Then 
\begin{align*}
    &\det L = \frac{1}{\det(L^{-1})} 
    = \prod \limits_{i=1}^r \frac{n^{|A_i|-1}}{|A_i|} 
    = \frac{n^{n-r}}{ \prod \limits_{i=1}^r |A_i|} 
    = \frac{n^{n-2}}{n ^{r-2} \prod \limits_{i=1}^r |A_i|}.
\end{align*}
By Theorem~\ref{main}, we have $\det L = {\sum_T w(T)}/{\sum_F w(F)},$
where the sum in the denominator is over all valid forests $F$ in the superport network, and the sum in the numerator is over all spanning trees $T$. By Corollary~\ref{Cayley}, $\sum_T w(T) = n^{n-2}$. Then $
n^ {r-2} \cdot \prod \limits_{i=1}^r |A_i| = {\sum_F w(F)}$, where the sum is over all valid forests $F$. Since the conductances are unit, this sum ${\sum_F w(F)}=\sum_F 1$ equals the number of valid forests. By the definition of a valid forest, this equals the number of desired trees.
\end{proof}

For the proof of Theorem~\ref{main}, we need several auxiliary notions and lemmas.

First we express $\det L$
in terms of so-called valid minors of the response matrix $C$ of the electrical network obtained by unifying all the superports.
A set $\{i_1, \dots, i_{p-1}\}$ of integers is \emph{valid} if $ i_k \in A_k$ for each $k=1,\dots,p-1$. A \emph{valid minor} is a minor $C_{I}^{J}$ formed by a valid set of rows $I$ and a valid set of columns $J$. \newremove{}
    

\begin{lemma} \label{det_L_in_terms_of_C} 
 Let $L$ and $C$ be the response matrices of a superport network \bluenew{$m > p$} and the electrical network obtained by unifying all the superports, respectively. 
 Then
	$$\det L= \frac{\det \widetilde{C}}{\sum_{I,J} \det C_I^J},$$
where the sum in the denominator is taken over all valid sets ${I}, {J}$.
\end{lemma}

\begin{proof} 
	Denote by \bluenew{$F$, $G$, and $H$}, respectively, the matrices that are obtained after steps 1, 2, and 3 of Algorithm~\ref{alg:C2L}. First, we write a chain of equalities that prove the lemma 
    and then explain each of them:	
\begin{align*}
		 \det L^{-1} & \stackrel{(1)}{=}
		\det H_{M \setminus R} ^ {M \setminus R} \\
        & \stackrel{(2)}{=} 
        \sum\limits_{ \bluevarnew{j_1} = 1}^{r_1} \sum\limits_{\bluevarnew{j_2} = r_1+1}^{r_2} \dots \sum\limits_{\bluevarnew{j_{p-1}} = r_{p-2}+1}^{r_{p-1}} 
        (-1)^{\sum\limits_{t=1}^{p-1}(r_t-\bluevarnew{j_t})} \det G_{M \setminus \bluevarnew{\{j_1, j_2, \dots, j_{p-1}\}}}^ {M \setminus R}
        \\
        & \stackrel{(3)}{=} 
        \sum\limits_{J \text{ valid}} (-1)^{\sum\limits_{t =1}^{p-1} r_t + \sum\limits_{\bluevarnew{j \in J}} \bluevarnew{j}} \det G_{M \setminus \bluevarnew{J}}^ {M \setminus R}
        \\
		& \stackrel{(4)}{=}
		\sum\limits_{I, J \text{ valid}} (-1)^{\sum\limits_{i \in I} i + \sum\limits_{j \in J} j} \det F_{M \setminus \bluevarnew{J}}^ {M \setminus \bluevarnew{I}} \\
		&	\stackrel{(5)}{=}
		\sum\limits_{I, J \text{ valid}} \frac{\det (F^{-1})_{I}^{J}}{\det F^{-1} } \\
            &   \stackrel{(6)}{=}
		 \frac{ \sum\limits_{I, J \text{ valid}} \det C_{I}^{J}} {\det \widetilde{C}}.
	\end{align*}
 	Let us \bluenew{explain the equalities}:
	\begin{enumerate}
		\item[(1)] This holds by steps 4 and 5 of Algorithm~\ref{alg:C2L}. Recall that $M$ and $R=\{r_1,\dots, r_p\}$ denote the set of boundary vertices and the set of roots respectively.
        \item[(2)] \bluenew{Let us restrict ourselves to the rectangular submatrix $G_{M \setminus R}^{M}$ of $G$, and see how it is modified in step 3 of the algorithm. Denote  by $g_1, \dots, g_{m}$ the rows of the submatrix 
        and by $e_1, \dots, e_{m-p}$ the standard basis in $\mathds{R}^{m-p}$. In step 3, we first subtract 
        $g_{r_1}$ from 
        $g_1, \dots, g_{r_1-1}$. By the multilinear expansion, this changes the minor $\det G_{M \setminus R}^{M \setminus R}$ to}
        %
        %
        \begin{align*}
        &\frac{(g_1 - g_{r_1}) \wedge (g_2 - g_{r_1}) \wedge \dots \wedge (g_{r_1-1} - g_{r_1}) \wedge (g_{r_1+1} \wedge \dots \wedge g_{m-1})}{e_1 \wedge e_2 \wedge \dots \wedge e_{m-p}} \\
        & {=}  \sum\limits_{j_1 =1}^{r_1}(-1)^{r_1-j_1}
        \frac{(g_1 \wedge  \dots \wedge g_{j_1-1} \wedge g_{j_1+1}\wedge  \dots \wedge g_{r_1}) \wedge (g_{r_1+1} \wedge \dots \wedge g_{m-1})}{e_1 \wedge e_2 \wedge \dots \wedge e_{m-p}} \\
        &{=} \sum\limits_{j_1 =1}^{r_1}(-1)^{r_1-j_1} \det G^{M \setminus R}_{M \setminus R  \cup \{r_1\} \setminus \{j_1\}},
        \end{align*}
        where $g_1 \wedge \dots \wedge g_{j_1-1} \wedge g_{j_1+1}\wedge \dots \wedge g_{r_1}$ denotes the exterior product of all rows from $g_1$ to $g_{ r_1}$ except for $g_{j_1}$ and $g_{r_1+1} \wedge \dots \wedge g_{m-1}$ denotes the exterior product of rows \bluenew{$g_k$ for all $k>r_1$ such that $k \notin R$}. 
        Similarly, if the row $g_{r_k}$ is subtracted from the rows
        $g_{r_{k-1}+1}, \dots, g_{r_k-1}$ for each $k = 1, \dots, p-1$, then the \bluenew{minor in question becomes} $$\det H_{M \setminus R} ^ {M \setminus R} = \sum\limits_{j_1 = 1}^{r_1} \sum\limits_{j_2 = r_1+1}^{r_2} \dots \sum\limits_{j_{p-1} = r_{p-2}+1}^{r_{p-1}} \prod\limits_{t=1}^{p-1}(-1)^{r_t-j_t} \det G_{M \setminus \{j_1, j_2, \dots, j_{p-1}\}}^ {M \setminus R}.$$
        \item[(3)] This follows by the substitution $J = \{j_1, j_2, \dots, j_{p-1}\}.$
        \item[(4)] This is obtained in the same way as equations (2)--(3).
        \item[(5)] This is the Jacobi lemma on the complementary minors of two reciprocal matrices.
        \item[(6)] This is true by step 1 of the algorithm. By Theorem~\ref{C-L}, the lemma follows.
        \end{enumerate}
        \vspace{-0.8cm}
        \end{proof}
        
Now we apply Kenyon--Wilson formula~\eqref{equation_Kenyon_Wilson} to compute valid minors. 

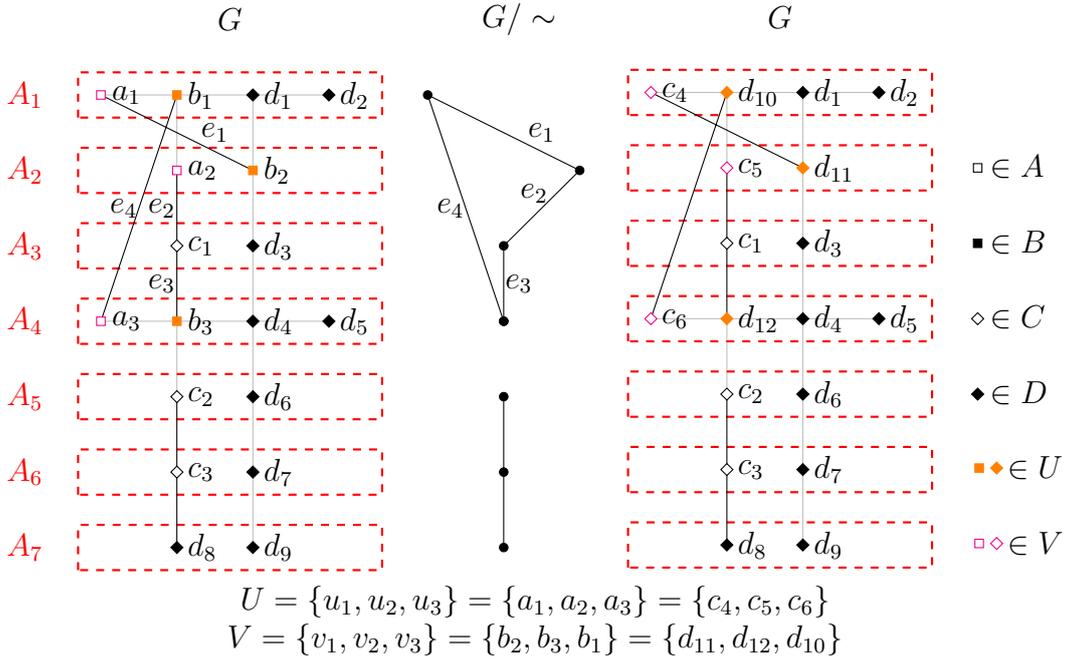
\begin{figure}[h]
    \centering
\begin{tabular}{c}
             \begin{circuitikz} 
             \draw[color=lightgray] (1, 2) to (1, 3)
             (1,5) to (1,6);
              \draw[color=lightgray] (2, 0) to (2, 1) to (2, 2) to (2, 3) to (2, 4) to (2, 5) to (2,6);
              \draw[color=lightgray] (0,3) to (1,3) to (2, 3) to (3,3);
              \draw[color=lightgray] (0,6) to (1,6) to (2, 6) to (3,6);
\draw[dashed, red, thick] (-0.3,-0.3) rectangle (3.7,0.3); 
\draw[dashed, red, thick] (-0.3,0.7) rectangle (3.7,1.3);
\draw[dashed, red, thick] (-0.3,1.7) rectangle (3.7,2.3);
\draw[dashed, red, thick] (-0.3,2.7) rectangle (3.7,3.3);
\draw[dashed, red, thick] (-0.3,3.7) rectangle (3.7,4.3);
\draw[dashed, red, thick] (-0.3,4.7) rectangle (3.7,5.3);
\draw[dashed, red, thick] (-0.3,5.7) rectangle (3.7,6.3);
             \draw
             (1,0) to (1,1) to (1,2)
             (0,3) to (1, 6)
             (1, 5) to (1, 4) to (1, 3)
             (0,6) to (2, 5)
             (1,0) node[diamondpole]{}
              (2,0) node[diamondpole]{}
              (1,1) node[odiamondpole]{}
              (2,1) node[diamondpole]{}
               (1,2) node[odiamondpole]{}
              (2,2) node[diamondpole]{}
               (1,3) node[squarepole, color=orange]{}
              (2,3) node[diamondpole]{}
              (0,3) node[osquarepole, color=magenta]{}
               (1,4) node[odiamondpole]{}
              (2,4) node[diamondpole]{}
              (3,3) node[diamondpole]{}
               (1,5) node[osquarepole, color=magenta]{}
              (2,5) node[squarepole, color=orange]{}
              (0,6) node[osquarepole, color=magenta]{}
              (1,6) node[squarepole, color=orange]{}
              (2,6) node[diamondpole]{}
              (3,6) node[diamondpole]{};
              \draw (-1,0) node[red]{$A_7$};
               \draw (-1,1) node[red]{$A_6$};
                \draw (-1,2) node[red]{$A_5$};
                 \draw (-1,3) node[red]{$A_4$};
                  \draw (-1,4) node[red]{$A_3$};
                   \draw (-1,5) node[red]{$A_2$};
                    \draw (-1,6) node[red]{$A_1$};
            \draw {[anchor= west]
            
                            (1,0) node[]{$w_8$}
                            (2,0) node []{$w_9$}
                             (1,1) node[]{$z_3$}
                            (2,1) node []{$w_7$}
                            (1,2) node[]{$z_2$}
                            (2,2) node []{$w_6$}
                            (0,3) node[]{$x_3$}
                            (1,3) node[]{$y_3$}
                            (2,3) node []{$w_4$}
                            (3,3) node[]{$w_5$}
                            (1,4) node[]{$z_1$}
                            (2,4) node []{$w_3$}
                            (1,5) node[]{$x_2$}
                            (2,5) node []{$y_2$}
                            (0,6) node[]{$x_1$}
                            (1,6) node[]{$y_1$}
                            (2,6) node []{$w_1$}
                            (3,6) node[]{$w_2$}
                            };
                            \draw (5.3,0) to (5.3,1) to (5.3,2)
                            (5.3, 3) to (5.3, 4) to (6.3, 5) to(4.3,6) to (5.3,3);
                            \draw (1.5, 5.5) node[]{$e_1$};
                            \draw (0.8, 4.5) node[]{$e_2$};
                            \draw (0.8, 3.5) node[]{$e_3$};
                            \draw (0.3, 4.5) node[]{$e_4$};

                            \draw (5.8, 5.5) node[]{$e_1$};
                            \draw (5.7, 4.7) node[]{$e_2$};
                            \draw (5.5, 3.5) node[]{$e_3$};
                            \draw (4.6, 4.5) node[]{$e_4$};
                            \draw (5.3,0) node[circ]{}
                            (5.3,1) node[circ]{}
                            (5.3,2) node[circ]{}
                            (5.3, 3) node[circ]{}
                            (5.3, 4) node[circ]{}
                            (6.3, 5) node[circ]{}
                            (4.3,6) node[circ]{}
                            (5.3,3) node[circ]{};
                            \draw (1.7, 7) node[]{$G$};
                            \draw (5.5, 7) node[]{$G/{\sim}$};
             \end{circuitikz}
\end{tabular}
\begin{tabular}{c}
             \begin{circuitikz} 
             \draw[color=lightgray] (1, 2) to (1, 3)
             (1,5) to (1,6);
              \draw[color=lightgray] (2, 0) to (2, 1) to (2, 2) to (2, 3) to (2, 4) to (2, 5) to (2,6);
              \draw[color=lightgray] (0,3) to (1,3) to (2, 3) to (3,3);
              \draw[color=lightgray] (0,6) to (1,6) to (2, 6) to (3,6);
\draw[dashed, red, thick] (-0.3,-0.3) rectangle (3.7,0.3); 
\draw[dashed, red, thick] (-0.3,0.7) rectangle (3.7,1.3);
\draw[dashed, red, thick] (-0.3,1.7) rectangle (3.7,2.3);
\draw[dashed, red, thick] (-0.3,2.7) rectangle (3.7,3.3);
\draw[dashed, red, thick] (-0.3,3.7) rectangle (3.7,4.3);
\draw[dashed, red, thick] (-0.3,4.7) rectangle (3.7,5.3);
\draw[dashed, red, thick] (-0.3,5.7) rectangle (3.7,6.3);
             \draw
             (1,0) to (1,1) to (1,2)
             (0,3) to (1, 6)
             (1, 5) to (1, 4) to (1, 3)
             (0,6) to (2, 5)
             (1,0) node[diamondpole]{}
              (2,0) node[diamondpole]{}
              (1,1) node[odiamondpole]{}
              (2,1) node[diamondpole]{}
               (1,2) node[odiamondpole]{}
              (2,2) node[diamondpole]{}
               (1,3) node[diamondpole, color=orange]{}
              (2,3) node[diamondpole]{}
              (0,3) node[odiamondpole, color=magenta]{}
               (1,4) node[odiamondpole]{}
              (2,4) node[diamondpole]{}
              (3,3) node[diamondpole]{}
               (1,5) node[odiamondpole, color=magenta]{}
              (2,5) node[diamondpole, color=orange]{}
              (0,6) node[odiamondpole, color=magenta]{}
              (1,6) node[diamondpole, color=orange]{}
              (2,6) node[diamondpole]{}
              (3,6) node[diamondpole]{};
            \draw {[anchor= west]
            
                            (1,0) node[]{$w_8$}
                            (2,0) node []{$w_9$}
                             (1,1) node[]{$z_3$}
                            (2,1) node []{$w_7$}
                            (1,2) node[]{$z_2$}
                            (2,2) node []{$w_6$}
                            (0,3) node[]{$z_6$}
                            (1,3) node[]{$w_{12}$}
                            (2,3) node []{$w_4$}
                            (3,3) node[]{$w_5$}
                            (1,4) node[]{$z_1$}
                            (2,4) node []{$w_3$}
                            (1,5) node[]{$z_5$}
                            (2,5) node []{$w_{11}$}
                            (0,6) node[]{$z_4$}
                            (1,6) node[]{$w_{10}$}
                            (2,6) node []{$w_1$}
                            (3,6) node[]{$w_2$}
                            };
                            \draw (4.3,2) node[diamondpole]{};
              \draw (4.3,3) node[odiamondpole]{};
              \draw (4.3,4) node[squarepole]{};
              \draw (4.3,5) node[osquarepole]{};
              \draw {[anchor= west]
              (4.32,5.02) node[]{ $\in X$}
              (4.32,4.02) node[]{ $\in Y$}
              (4.32,2.02) node[]{ $\in W$}
              (4.32,3.02) node[]{ $\in Z$}
              (4.55,1.02) node[]{ $\in U$}
              (4.55,0.02) node[]{ $\in V$}
              };
              \draw (4.55,1.02) node[diamondpole, color=orange]{};
              \draw (4.32,1.02) node[squarepole, color=orange]{};
              \draw (4.55,0.02) node[odiamondpole, color=magenta]{};
              \draw (4.32,0.02) node[osquarepole, color=magenta]{};
              \draw (1.7, 7) node[]{$G$};
             \end{circuitikz}
\end{tabular} 
$U=\{ u_1, u_2, u_3\}=\{x_1, x_2, x_3\}=\{z_4, z_5, z_6\}$\\
$V=\{ v_1, v_2, v_3\}=\{y_2, y_3, y_1\}=\{w_{11}, w_{12}, w_{10}\}$
\caption{Left: A graphical description of all the valid minors which a given forest \bluevarnew{contributes} to. The network, the superports $A_1,\dots,A_7$, and the forest $G$ are shown in gray, red, and black respectively. The boundary vertices are decomposed (coloured) into four sets (colours) $X,Y,Z,W$, depicted by \bluevarnew{empty squares, filled squares, empty rhombi, filled rhombi} respectively. See Lemma~\ref{equivalence_lemma}. Middle: The quotient $G/{\sim}$ contains a cycle with the edges $e_1,e_2,e_3,e_4$. See Lemma~\ref{cycle_lemma}. Right: A new decomposition (colouring) obtained by the involution~$f$. The vertices in the sets $U$ and $V$, shown in \bluevarnew{orange and pink} respectively, have changed their sets (colours). See case 1 of the proof of Lemma~\ref{main_lemma}. 
}
\label{fig:case1}
\end{figure}

\bluenew{Notice that for any valid sets 
$I, J$ and $X := I \setminus J$, $Y := J \setminus I$, $Z := I \cap J$, $W := M \setminus (I \cup J)$  
we have $\det C^J_I=\det C^{Y \cup Z}_{X \cup Z}=\det C^{Y, Z}_{X,Z}$.}
We say that a forest \emph{contributes} to \bluenew{the 
minor $C^J_I$} 
if this forest \bluenew{participates in some term in the sum} in the numerator of the right side of~\eqref{equation_Kenyon_Wilson} from Theorem~\ref{Kenyon-Wilson} (recall that $w(\dots|\dots|\dots)$ is defined as a sum over forests). The \emph{sign of the forest in the minor} is the sign with which the weight of this forest is included in this sum. In other words, the sign of the forest is $(-1)^{|X|} \cdot \mathrm{sgn}(\pi),$ where the set $X$ and the permutation $\pi$ are determined by the minor and how the forest connects the boundary vertices. We introduce a graphical description of all the valid minors which a given forest contributes to, in terms of the coloring of boundary vertices into four colors $X, Y, Z, W$. See Fig.~\ref{fig:case1} to the left.

\begin{lemma} \label{equivalence_lemma}
Let $G = G_1 \sqcup \dots \sqcup G_{c}$ be a spanning forest with the connected components $G_1, \dots, G_{c}$.
Then the valid minors which the forest $G$ contributes to are in a bijection with the partitions of the set $M$ into four subsets
$X, Y, Z, W$ satisfying the following conditions:
\begin{enumerate}
		\item[(1)] $A_p \subset W$;
		\item[(2)] for each $k= 1, \dots, p-1$ we have $|X \cap A_k| = |Y \cap A_k| = 1-|Z \cap A_k| = 0$ or $1$;
		\item[(3)] for each $k= 1, \dots, c$ we have $|X \cap G_k| = |Y \cap G_k| = 1-|W \cap G_k| = 0$ or $1$.
\end{enumerate}

Moreover, the sign of the forest $G$ in the valid minor is 
\begin{equation}\label{eq-sgn-abcdg}
\mathrm{sgn}(G, X, Y, Z, W) :=(-1)^{|X|} \mathrm{sgn}(\sigma \circ \tau),
\end{equation}
where the bijections $\sigma\colon X \to Y$ and $\tau\colon Y \to X$ are uniquely determined by the following two conditions:
\begin{enumerate}
    \item[($\bluevarnew{\sigma}$)] $x$ and $\sigma(x)$ belong to the same connected component \bluenew{of} $G$ for each $x \in X$;
    \item[($\bluevarnew{\tau}$)] $y$ and $\tau(y)$ belong to the same superport for each $y \in Y$.
\end{enumerate}
\end{lemma}

\bluenew{\begin{remark} \label{remark_about_c}
From conditions (1)--(3), we can deduce that $c = m - p + 1$. 
Indeed, according to condition (3), each component contains exactly one vertex from either $W$ or $X$, which gives $c = |W| + |X|$.
Conditions (1) and (2) imply that all superports, except for the last one, contain exactly one vertex from either $Y$ or $Z$. Therefore, we have $|Y| + |Z| = p - 1$.
Combining these results, we find
$c = |W| + |X| = |M| - |Y| - |Z| = m - p + 1.$
\end{remark}}

Recall that  the sign of an empty permutation is 
$+1$; this happens when $X = Y=\emptyset$.

\begin{proof}
\bluenew{First, we explain how the bijections $\sigma$ and $\tau$ are defined and why they exist. For each $x \in X$, the element $\sigma(x) \in Y$ is the unique element in the same connected component of $G$ as $x$. By condition (3), in each connected component $G_k$, if there is a vertex from $X$ or $Y$, then exactly one vertex belongs to $X$ and exactly one to $Y$. This ensures that $\sigma$ is a well-defined bijection. Similarly, $\tau$ is a well-defined bijection by conditions (1)--(2).}

Now, we prove that if $G$ contributes to a valid minor $C_{I}^{J},$ then conditions (1)--(3) are satisfied for $Z := I \cap J$, $X := I \setminus Z$, $Y := J \setminus Z$, and $W := M \setminus (X \cup Y \cup Z)$.

Since the minor $C_{I}^{J}$ is valid, it follows that $I = \{i_1, \dots, i_{p-1}\}$ and $ J= \{j_1, \dots, j_{p- 1}\}$ satisfy the condition $i_k, j_k \in A_k$ for all $k=1,\dots,p-1$. Then $A_p \cap (I \cup J)=\emptyset$. By our construction, $X \cup Y \cup Z = I \cup J$, hence $A_p \subset W$ and condition (1) holds.

By our construction, $(X \cup Y \cup Z) \cap A_k = \{i_k, j_k\}$. If $i_k = j_k$, then $(X \cup Y \cup Z) \cap A_k = \{ i_k \} \subset Z$ and $|X \cap A_k|$ = $|Y \cap A_k|$ = $1-|Z \cap A_k| = 0$. If $i_k \neq j_k$, then $i_k \in X$ and $j_k \in Y$, and $|X \cap A_k|$ = $|Y \cap A_k|$ = $1-|Z \cap A_k| = 1.$ In both cases, condition (2) holds.

By Theorem~\ref{Kenyon-Wilson}, condition (3) is satisfied and the forest contributes to the minor with the sign $(-1)^{|X|} \mathrm{sgn}(\sigma \circ \tau)$ because \bluenew{$\rho^{-1} \circ \sigma \circ \tau \circ \rho = \pi$, where we set $X = \{x_1, \dots, x_{|X|}\}$ with $x_1 < \dots < x_{|X|}$, $Y = \{y_1, \dots, y_{|Y|}\}$ with $y_1 < \dots < y_{|Y|}$, and $\rho(i) := x_i$}.

\bluenew{Finally,} we prove that if conditions (1)--(3) are satisfied, then $G$ contributes to the minor $C_{I}^{J},$ where $I:=X \cup Z, J:=Y \cup Z$, and this minor is valid. It follows from conditions (1) and (2) that the minor is valid. From condition (3), by Theorem~\ref{Kenyon-Wilson} we get that $G$ contributes to the minor with the sign $(-1)^{|X|} \mathrm{sgn}(\sigma \circ \tau)$.
\end{proof}

\begin{corollary} \label{equivalence_corollary}
    For any superport network, we have $$\sum\limits_{I, J \text{ valid}} \det C_I^J= \frac{ \sum \limits_{X, Y, Z, W, G} \mathrm{sgn}(G, X, Y, Z, W) \cdot w(G)}{\sum\limits_H w(H)},$$ where the sum in the denominator is over all valid forests $H$ in the electric network obtained by unifying all the superports, the sum in the numerator is over the partitions $X \sqcup Y \sqcup Z \sqcup W$ of $M$ and spanning forests $G$ satisfying conditions (1)--(3) of Lemma~\ref{equivalence_lemma}, and we use notation~\eqref{eq-sgn-abcdg}.
\end{corollary}

This corollary allows working with partitions $M = X \sqcup Y \sqcup Z \sqcup W$ instead of minors. For non-valid forests $G$, there is the following strong restriction on such partitions.

\begin{lemma} \label{cycle_lemma} \textup{(See Fig.~\ref{fig:case1})}
Let a spanning forest $G$ and a partition $X \sqcup Y \sqcup Z \sqcup W$ of $M$ satisfy conditions (1)--(3) of Lemma~\ref{equivalence_lemma}.
Let $e_1,...,e_n$ be a sequence of oriented edges forming a simple oriented cycle in the quotient $G/{\sim}$. In the forest $G$, the union $e_1 \cup \dots \cup e_n$ splits into disjoint paths denoted by $u_1 \dots v_1, \dots, u_l \dots v_l$ in the order they appear in the cycle. Then up to the reversal of the cycle direction one of the following conditions holds:
\begin{enumerate}
	\item[(XY)] $u_k\in X$ and $v_k\in Y$ for each $k=1, \dots, l$; or
	\item[(ZW)] $u_k\in Z$ and $v_k \in W$ for each $k=1, \dots, l$.
\end{enumerate}
\end{lemma}

\begin{proof}
For convenience, denote $v_{l+1}:=v_1$ and $u_{l+1}:=u_1$. For each $k=1, \dots, l$ the vertices $u_k$ and $v_k$ are in the same connected component ($u_k, v_k \in G_i$ for some $i$) and $v_k$ and $u_{k+1}$ are in the same superport ($v_k, u_{k+1} \in A_j$ for some $j$).

Let us prove that if some $v_k \in Z,$ then $u_{k+1} \in W$ and $v_{k+1}\in Z.$ Indeed, since $v_k \in Z$ and $v_k, u_{k+1} \in A_j$, by condition (1) it follows that $j \neq p$ and by condition (2) it follows that $u_{k+1} \in W$. Since $u_{k+1} \in W$ and $u_{k+1}, v_{k+1} \in G_i,$ by condition (3) it follows that $v_{k+1}\in Z.$ Thus, if some $v_k \in Z,$ then condition (ZW) holds \bluenew{after changing the direction of the cycle.}

Similarly, if some $v_k \in W$, then $u_k \in Z$ and $v_{k-1} \in W$. Then condition (ZW) \bluenew{holds}.

It remains to consider the case when all $u_{k}, v_{k} \notin Z \cup W$, i.e. all $u_{k}, v_{k} \in X \cup Y$. Then by (1)--(3) it follows that condition (XY) is fulfilled up to change of the cycle direction.
\end{proof}

We now come to the key step of the proof: cancellation of non-valid forests.

\begin{lemma} \label{main_lemma} We have
\begin{align}\label{signed_formula}
    \sum_{X, Y, Z, W, G} \mathrm{sgn}(G, X, Y, Z, W) \cdot w(G) = \sum_F w(F),
\end{align}
 where the first sum is over spanning forests $G$ and the partitions $X \sqcup Y \sqcup Z \sqcup W$ of $M$ satisfying conditions (1)--(3) of Lemma~\ref{equivalence_lemma}, the second sum is over all valid forests $F$, and we use notation~\eqref{eq-sgn-abcdg}.
\end{lemma}

\begin{proof}
Let us show that the contributions of all non-valid forests to the left side \bluenew{of (\ref{signed_formula})} cancel, and valid forests contribute exactly once with a positive sign.

\bluenew{By Remark~\ref{remark_about_c}, $c = m-p+1$ for any spanning forest $G=G_1\sqcup\dots\sqcup G_c$ that satisfies conditions (1)--(3).} Fix a spanning forest $G$ with $m-p+1$ connected components. Such a forest $G$ is valid if and only if $G/{\sim}$ has no cycles. Indeed, since $G$ has exactly $m-p+1$ connected components and $n$ vertices, it has $p-1 + n-m$ edges. Then $G/{\sim}$ also has $p-1 + n-m$ edges and $p + n-m$ vertices. Thus either $G/{\sim}$ is a spanning tree and $G$ is a valid forest, or $G/{\sim}$ has a cycle.
Consider these two cases separately. Consider the \bluenew{bijections} $\sigma$ and $\tau$ defined in Lemma~\ref{equivalence_lemma}.


\emph{Case 1}: $G$ is not valid, i.e., \bluenew{$G/{\sim}$} has a cycle. See Fig.~\ref{fig:case1}. Let us show that $G$ contributes the same number of times with the signs $+1$ and $-1$ to the left-hand side of~\eqref{signed_formula}. We construct a sign-reversing involution on the set of partitions $X\sqcup Y\sqcup Z \sqcup W = M$ satisfying conditions (1)--(3).

For the construction, fix a canonical choice of a cycle in $G/{\sim}$ as follows. Since the vertices in our network are \bluenew{numbered}, we can write each oriented cycle as a string of vertex numbers, starting from the vertex with the smallest number. 
The simple oriented cycle in $G/{\sim}$ that corresponds to the lexicographically smallest string is called \emph{the main cycle}.

The union of edges forming the main cycle in $G/{\sim}$ splits into disjoint simple paths $u_1 \dots v_1, \dots, u_l \dots v_l$ in $G$. Set $U := \{u_1, \dots, u_l\}$, and $V := \{v_1, \dots, v_l\}$. By Lemma~\ref{cycle_lemma}, changing the cycle direction, if necessary, we may assume that one of two conditions (XY) or (ZW) is satisfied, i.e. either $U \subset X, V \subset Y$ or $U \subset Z, V \subset W$.

Consider the following involution $f$ on the set of partitions: 
$$
	f(X, Y, Z, W) =
	\begin{cases}
		(X \setminus U, Y \setminus V, Z \cup U, W \cup V), &\text{if $U \subset X, V \subset Y$};\\
		(X \cup U , Y \cup V, Z \setminus U, W \setminus V), &\text{if $U \subset Z, V \subset W$}.
	\end{cases}
	$$

First let us show that $f$ preserves conditions (1)--(3).
By condition (1) and Lemma~\ref{cycle_lemma}, $A_p$ has no vertices in $U\cup V$, hence condition (1) is preserved.
By conditions (2), (3), and Lemma~\ref{cycle_lemma} each superport (except the last one) and each connected component of $G$ contains either one vertex in $U$ and one in $V$, or no vertices in $U\cup V$.
If $k$-th superport or connected component has no vertices in $U\cup V$, then conditions (2) and (3) for this $k$ are preserved automatically. If $k$-th superport or connected component has one vertex in $U$ and one in $V$, then in conditions (2) and (3) for this $k$ equality to zero changes to equality to one and vice versa.

Now check the identity $f(f(X, Y, Z, W)) = (X, Y, Z, W)$. Applying the mapping $f$ twice, we get
\begin{align*}
	&(X, Y, Z, W) \mapsto
	\begin{cases}
		(X \setminus U, Y \setminus V, Z \cup U, W \cup V), &\text{if $U \subset X, V \subset Y$};\\
		(X \cup U , Y \cup V, Z \setminus U, W \setminus V), &\text{if $U \subset Z, V \subset W$};
	\end{cases}
 \\ &\mapsto \begin{cases}
		((X \setminus U) \cup U, (Y \setminus V) \cup V, (Z \cup U) \setminus U, (W \cup V) \setminus V), &\text{if $U \subset X, V \subset Y$};\\
		((X \cup U) \setminus U, (Y \cup V) \setminus V, (Z \setminus U) \cup U, (W \setminus V) \cup V), &\text{if $U \subset Z, V \subset W$};
	\end{cases} = (X, Y, Z, W).
\end{align*}
Thus, $f$ is an involution.

It remains to show that $f$ reverses the sign of the forest $G$. By Lemma~\ref{equivalence_lemma}, the sign is $(-1)^{|X|} \mathrm{sgn}(\sigma \circ \tau).$ After mapping $f$, the number $|X|$ is replaced by $|X| \pm |U|$ and the cycle $(u_1\dots u_{|U|})$ of length \bluenew{$|U| \neq 0$} is removed from or added to the permutation $\sigma \circ \tau$. Thus the permutation sign is multiplied by $(-1)^{|U|+1}$ because an even-length cycle is an odd permutation and an odd-length cycle is an even one. We get 
\begin{multline*}\mathrm{sgn}(G, f(X, Y, Z, W))=(-1)^{|X| \pm |U|} \mathrm{sgn}(\sigma \circ \tau) \cdot (-1)^{|U|+1} =\\= (-1)^{|X| + 1} \mathrm{sgn}(\sigma \circ \tau) = -\mathrm{sgn}(G, X, Y, Z, W).
\end{multline*}

\emph{Case 2}: $G$ is valid, i.e. there is no cycle in $G/{\sim}$. See Fig.~\ref{fig:last}. Let us show that $G$ contributes to the left-hand side of~\eqref{signed_formula} with the coefficient +1; that is, the partition $X \sqcup Y \sqcup Z \sqcup W$ is uniquely determined by $G$ and conditions (1)--(3), and we have $\mathrm{sgn}(G, X, Y, Z, W)=+1$.

Let us prove that $X = \bluevarnew{\emptyset}$. Assume the converse. Take any element $x_1$ from $X$. Set $y_i = \sigma(x_i)$ for all $i$ and $x_i = \tau(y_{i-1})$ for $i > 1.$ 
We get a path \bluenew{$[x_1] [y_1] [x_2] [y_2] \dots$} in $G/{\sim}$. Since there is a finite number of vertices, it follows that $x_i = x_1$ for some \bluenew{$i\ne 1$. Take the minimal $i$ with this property. Then $[x_1] [y_1] \dots [y_{i-1}]$ is a simple cycle in $G/{\sim}$ because the paths $x_j y_j$ are disjoint in $G$ by condition (3).} This contradiction proves that $X = \bluevarnew{\emptyset}$.

In the tree $G/{\sim},$ mark the vertex corresponding to the last superport. Let it be the root of the tree. Introduce a strict partial order on vertices of the graph $G$ letting $u \prec v$ if in the tree $G/{\sim}$, vertices $[u]$ and $[v]$ are different and $[u]$ belongs to the unique \bluenew{simple} path connecting the root with the vertex $[v]$.

Since $X =  \bluevarnew{\emptyset}$, by condition (3), we get $Y = \bluevarnew{\emptyset}$ and for each $k$ we have $|W \cap G_k| = 1$. That is, each component of $G$ has exactly one vertex in $W$ and all the other vertices belong to $Z$.

\begin{figure}[h]
\center
\begin{tabular}{c}
             \begin{circuitikz} 
             \draw[color=lightgray] (0,0) to (0, 1)  (0,2) to (0,3) to (0,4) to (0,5) to (0,6);
             \draw[color=lightgray] (1,1) to (1,2) to (1,3) to (1,4) to (1,5) to (1,6);
\draw[dashed, red, thick] (-0.5,-0.3) rectangle (1.5,0.3); 
\draw[dashed, red, thick] (-0.5,0.7) rectangle (1.5,1.3);
\draw[dashed, red, thick] (-0.5,1.7) rectangle (1.5,2.3);
\draw[dashed, red, thick] (-0.5,2.7) rectangle (1.5,3.3);
\draw[dashed, red, thick] (-0.5,3.7) rectangle (1.5,4.3);
\draw[dashed, red, thick] (-0.5,4.7) rectangle (1.5,5.3);
\draw[dashed, red, thick] (-0.5,5.7) rectangle (1.5,6.3);
             \draw
              (0,6) to (1,4)
              (0,5) to (1,4) 
              (1,4) to (1, 3)
              (0,3) to[short, l_=${e_3}$] (0,2) to[short, l_=${e_2}$] (0,1)
              (1, 1) to[short, l_=${e_1}$] (1, 0);
              \draw
              (0,0) node[diamondpole]{}
              (0,1) node[diamondpole]{}
              (0,2) node[odiamondpole]{}
              (0,3) node[odiamondpole]{}
              (0,4) node[diamondpole]{}
              (0,5) node[odiamondpole]{}
              (0,6) node[odiamondpole]{}
              (1,0) node[diamondpole]{}
              (1,1) node[odiamondpole]{}
              (1,2) node[diamondpole]{}
              (1,3) node[diamondpole, color=red]{}
              (1,4) node[odiamondpole]{}
              (1,5) node[diamondpole]{}
              (1,6) node[diamondpole]{};
              \draw (-1,0) node[red]{$A_7$};
               \draw (-1,1) node[red]{$A_6$};
                \draw (-1,2) node[red]{$A_5$};
                 \draw (-1,3) node[red]{$A_4$};
                  \draw (-1,4) node[red]{$A_3$};
                   \draw (-1,5) node[red]{$A_2$};
                    \draw (-1,6) node[red]{$A_1$};
            \draw {[anchor= west]
              (0,1) node[]{$q$}
              (0,2) node[]{$u$}
              (0,3) node[]{$v$}
              (1,3) node[color=red]{$w$}
              };

              \draw (2,3) node[diamondpole]{};
              \draw (2,4) node[odiamondpole]{};
              \draw {[anchor= west]
              (2.02,3.02) node[]{ $\in W$}
              (2.02,4.02) node[]{ $\in Z$}};
              \draw (0.5, 6.6) node[]{$G$};
             \end{circuitikz}
\end{tabular} 
 \begin{tabular}{l}
             \begin{circuitikz} 
             
              \draw (1,0) to[short, l=${e_1}$] (1,1) to[short, l=${e_2}$] (1,2) to[short, l=${e_3}$] (1, 3) to (1,4)
              (1,4) to (0,5)
              (1,4) to (0,6) ;
              \draw
              (1,0)node[circ]{}
              (1,1)node[circ]{}
              (1,2)node[circ]{}
              (1,3)node[circ, color=red]{}
              (1,4)node[circ]{}
              (0,5)node[circ]{}
              (0,6)node[circ]{};
              \draw {[anchor= west]
              (1,3) node[color=red]{$[w]=[v]$}
              (1,2) node[]{$[u]$}
              (1,1) node[]{$[q]$}
              };
              \draw (0, -0.3) node[]{};
            \draw (0.5, 6.6) node[]{$G/{\sim}$};
             \end{circuitikz}
\end{tabular} 
\caption{A valid forest $G$ (left) and the tree $G/{\sim}$ (right). \bluevarnew{The empty and filled rhombi} depict the vertices of $Z$ and $W$ respectively. The vertices of the tree $G/{\sim}$ are partially ordered from bottom to top. 
See Case 2 of the proof of Lemma~\ref{main_lemma}. 
}
\label{fig:last}
\end{figure}
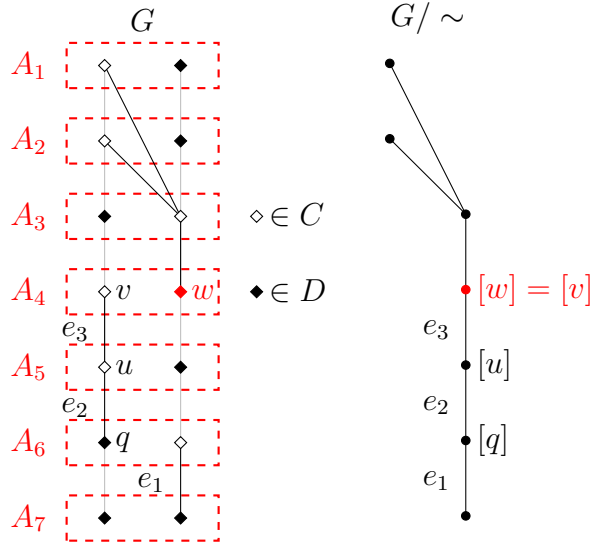

Let us prove that the least element in each component belongs to $W$ and all the others belong to $Z$. We illustrate the argument by an example in Fig.~\ref{fig:last}. Assume the converse. 
Then there is an element that is the least in its component in $G$ and is not in $W$. Among all such elements, we take some minimal element 
$w$. By our assumption, we get $w \in Z$ and by condition (1) we get 
$w \notin A_p$. Then there exists a \bluenew{unique nonempty path} connecting the root with the vertex $[w]$ in the tree $G/{\sim}$. Let $e_1, \dots, e_n$ be the \bluenew{consecutive} edges of the path.
Consider the same edges in the forest $G$. Take the component containing the last edge $e_n$. Denote by $q$ the least element of this component. The edge $e_n$ connects some vertices $u \prec v$, where $[v] = [w]$. Thus $v \neq w$ and $v \neq q$, because $w$ and $q$ are the least elements in their components. By condition (2) we get $v \in W$ because $w \in Z$ and vertices $v$ and $w$ are in the same superport. By condition (3) we get $q \in Z$ because $v \in W$ and vertices $q$ and $v$ are in the same component. Moreover, $q \prec w$ because $q \prec v$ and $[v] = [w]$. This contradicts the minimality of $w$, and proves that the least element in each component of $G$ belongs to $W$ and all the others belong to $Z$.

We obtain that the set $W$, hence $Z$, is uniquely determined by the forest $G$. \bluenew{Clearly, the resulting decomposition $M = Z \sqcup W$ satisfies conditions (1)--(3).} Moreover, \bluevarnew{we have} $\mathrm{sgn}(G, X, Y, Z, W) = (-1)^{|X|} \mathrm{sgn}(\sigma \circ \tau) = 1,$ because $X = \bluevarnew{\emptyset}$ and the permutation is empty.

Thus each valid forest $G$ contributes to the left side of~\eqref{signed_formula} with the coefficient $+1$.
\end{proof}

\begin{proof}[Proof of Theorem~\ref{main}] 
Combining the identities of Lemmas~\ref{det_L_in_terms_of_C}, \ref{main_lemma}, Corollary~\ref{equivalence_corollary}, and part (1) of Theorem~\ref{Kirchhoff}, we get $$\det L= \frac{\det \widetilde{C}}{\sum\limits_{I, J \text{ valid}}\det C_I^J} =  \frac{\sum\limits_T w(T)}{\sum\limits_{X, Y, Z, W, G} \mathrm{sgn}(G, X, Y, Z, W)w(G)}  = \frac{\sum\limits_T w(T)}{\sum\limits_F w(F)},$$
where 
\begin{itemize}
\item the sum $\sum\limits_{X, Y, Z, W, G} \mathrm{sgn}(G, X, Y, Z, W)w(G)$ is over spanning forests $G$ and the partitions $X\sqcup Y\sqcup Z \sqcup W$ of $M$ satisfying conditions (1)--(3) of Lemma~\ref{equivalence_lemma};
\item the sum $\sum\limits_T w(T)$ is over all spanning trees $T$ in the superport network;
\item the sum $\sum\limits_F w(F)$ is over all valid forests $F$ in the superport network.
\end{itemize}
\vspace{-1.0cm}
\end{proof}





\notformoebius
\section{Open problems} \label{sec-open}

Superport networks form a natural general setup to develop network theory. Most of the natural questions in this setup are open.

The most shouting open problem is to generalize Theorems~\ref{main} and~\ref{th-Lij} to the minors of arbitrary size, or in other words, generalize Theorem~\ref{Kenyon-Wilson} to superport networks.

\begin{problem}[All-minors matrix-tree theorem] Express the minors (of an arbitrary size) of the response matrix of a superport network in terms of sums over spanning forests.
\end{problem}

The main difficulty here is to determine the signs with which the forests contribute. For discrete period matrices, this has been achieved recently \cite{Lam-etal-25}.

The inverse problem for superport networks is also widely open; 
cf.~\bluenew{\cite{CM, Milton-Seppecher-08, Skopenkov-15}}:

\begin{problem}[Inverse problem] Characterize possible responses of:
\begin{enumerate}
    \item\label{item-1} superport networks with given sizes of superports;
    \item planar superport networks (defined analogously to~\cite{CM});
    \item\label{item-3} superport networks with complex conductances having positive real parts (cf.~\cite{Rote}).
\end{enumerate}
For each of the classes~\ref{item-1}--\ref{item-3}, give a construction algorithm of a superport network with a given response matrix.
\end{problem}

We remark that complex conductances with positive real parts describe alternating-current networks, common in electrical engineering. See \cite{PS} for an introduction.

A related problem is to study network transformations. Superport networks have their own transformations, not present for ordinary electrical networks; see e.g. Figure~\ref{fig:Box-H_transformation}.
\newremove{}

\begin{problem}[Network transformations] \newremove{} Find a set of transformations preserving the response and sufficient to transform any 
superport network to any other 
\bluenew{one} with the same \bluevarnew{superports and} response.
\end{problem}

It is interesting to develop an analogous theory and find its applications to other models on graphs; this could reveal new connections between the models~\cite{Galashin-Pylyavskyy-20,KW06,Lam-18,Skopenkov-23,Skopenkov-22}:

\begin{problem}[Other models]
    What are the analogs of superport boundary conditions (P) and (B) for simple random walks? for the Ising model? for the dimer model? for lattice gauge theories?
\end{problem}

Superport networks have already been successfully applied to the characterization of square-tilable polygons \cite{S}. They are a promising tool for the following open problem.

\begin{problem}[Square-tilable surfaces]
    Characterize which surfaces with a polyhedral metric (possibly, with boundary) are square-tilable (cf.~\cite{Dmitriev-Ozhegov, K, Novikov}).
\end{problem}
\endnotformoebius

\subsection*{Acknowledgements}
Our sincere thanks go to Omer Angel, Richard Kenyon, Asaf Nachmias, \bluer{Ivan Novikov}, Fedor Ozhegov, and Alexey Ustinov for taking the time to listen about this research and providing valuable comments and questions that have helped refining this paper.

\bibliographystyle{elsarticle-num}

{\small

\noindent
\textsc{Pavlo Pylyavskyy\\
University of Minnesota}\\
\texttt{ppylyavs\,@\,umn$\cdot $edu} \quad \url{https://sites.google.com/site/pylyavskyy/home}

\medskip

\noindent
\textsc{Svetlana Shirokovskikh\\
HSE University}\\
\texttt{sveta.17.10\,@\,yandex$\cdot $ru}

\medskip
\newremove{}
\noindent
\textsc{Mikhail Skopenkov\\
King Abdullah University of Science and Technology}
\\
\texttt{mikhail.skopenkov\,@\,gmail$\cdot $com} \quad \url{https://users.mccme.ru/mskopenkov/}

}

\end{document}